\documentclass{article}
\usepackage[left=2cm,right=2cm,top=2cm,bottom=2cm, marginparwidth=4cm, marginparsep=0.5cm]{geometry} 
\usepackage{amsmath,amssymb,mathtools} 
\usepackage{amsthm}
\usepackage{graphicx}
\usepackage{mathrsfs}
\usepackage{marginnote} 
\usepackage[utf8]{inputenc}         
\usepackage[T1]{fontenc}            %
\usepackage[final]{microtype} 
\usepackage{url}                    
\usepackage{tikz-cd}                   
\usetikzlibrary{arrows}                 
\usepackage{adjustbox}
\usepackage[]{color} 
\usepackage{framed} 
\usepackage[hidelinks]{hyperref}
\hypersetup{
  colorlinks   = true, 
  urlcolor     = blue, 
  linkcolor    = blue, 
  citecolor   = blue 
}
\usepackage[capitalise]{cleveref}

\tikzset{
commutative diagrams/.cd,
arrow style=tikz,
diagrams={>=latex}} 

\usepackage{titlesec} 
\titleformat*{\section}{\centering\large\sc }
\titleformat{\subsection}[runin]{\bfseries}{\thesubsection.}{3pt}{}

\usepackage{enumitem} 

\usepackage{tocloft} 

\theoremstyle{definition}
\newtheorem{thm}{Theorem}[section]
\newtheorem{defn}[thm]{Definition}
\newtheorem{rke}[thm]{Remark}
\newtheorem{ex}[thm]{Example}
\newtheorem{prop}[thm]{Proposition}

\newtheorem{lm}[thm]{Lemma}
\newtheorem{corr}[thm]{Corollary}

\newtheorem{con}[thm]{Conjecture}
\newtheorem*{ththm}{Theorem}
\newtheorem*{thcon}{Conjecture}

\newtheoremstyle{break} 
  {\topsep}{\topsep}%
  {\itshape}{}%
  {\bfseries}{}%
  {\newline}{}%
\theoremstyle{break}

\newcommand{\m}{\mathfrak{m}}
\newcommand{\ie}{\emph{i.e.,} }
\newcommand{\cf}{\emph{cf.} }

\newcommand{\R}{\mathbb{R}}
\newcommand{\C}{\mathbb{C}}
\newcommand{\Q}{\mathbb{Q}}
\newcommand{\Z}{\mathbb{Z}}
\newcommand{\N}{\mathbb{N}}

\newcommand{\Or}{\mathscr{O}}

\newcommand{\init}{\text{in}}

\begin{document}

\title{Minimal Pairs, Truncations and Diskoids}
\author{Andrei Bengu\c s-Lasnier}
\date{}
\maketitle

\begin{abstract}
    We build on the correspondence between abstract key polynomials and minimal 
    pairs made by Novacoski and show how to relate the valuations that are
    generated by each object. We can then give a geometric interpretation of
    valuations built in this fashion. To do so we employ an object called diskoid,
    which is a generalisation of the classical concept of ball in non-archimedian
    valued fields.
\end{abstract}
\setcounter{tocdepth}{1}
\tableofcontents
\vspace{1cm}

\section*{Introduction}

Extending valuations from a field $K$ onto a simple polynomial ring $K[X]$ can be
encoded in different ways and following different strategies. MacLane first introduced
key polynomials and augmented valuations in order to study simple extensions of
discretely valued fields of rank $1$ \cite{McL1}. This allowed him to prove results
in classical valuation theory \cite{McL2}. As the problem of  
local uniformization began to raise renewed interest, so did the study of the
extension problem of valuations from $K$ to $K[X]$. Problems of ramifications and 
the study of defect of extensions of valuations are deeply related to local
uniformization as it has been shown in \cite{CP}, \cite{CM}.

The problem
became an inductive step in understanding the valuations of function fields of
algebraic varieties. In dimension $2$, the problem has a geometric flavour:

\begin{itemize}[noitemsep]
    \item[--] The theory of plane curves and jet schemes has been used to give a
    precise description of minimal generating sequences of valuations in dimension
    2 (\cite{Te1}, \cite{Ab1}, \cite{Mo}, \cite{Tree}).
    \item[--] Curvettes are being used in \cite{Sp} whose results have later been
    generalised to non-algebraically closed fields by Cutkosky and Vinh in
    \cite{CV}.
\end{itemize}

Our goal in this paper is to establish a correspondence between certain key 
polynomials and geometric objects that would give a faithful correspondence between
the valuations induced by these key polynomials and the valuations given by these
geometric objects. We briefly mention the various notions that are involved in our
study starting with key polynomials. Different notions of key polynomials appeared
since MacLane's first definition in 1936. Vaquié generalised this first approach in
a series of articles (\cite{Va1}, \cite{Va2}, \cite{Va3}) where he defines limit
key polynomial and limit valuation in order to extend MacLane's theory to
non-necessarily discrete valuations of rank $1$. Independently Spivakovsky, Olalla
Acosta and Hererra Govantes gave a different approach to key polynomials, in
\cite{HOS},\cite{HOMS}. 
An alternative approach of key polynomials, which will be of interest to us, was
introduced in \cite{DMS}, that  of abstract key polynomials (see Section 1).
These allow to construct valuations, via a process called truncation: from $\mu$
a valuation on $K[X]$ and a key polynomial $Q$, we consider for any $f\in K[X]$, 
its $Q$-expansion
\[f=f_0+f_1 Q+\ldots+f_n Q^n,\;\forall i,\,\deg f_i<\deg Q\]
and define a map $\mu_Q$, by setting
\[\mu_Q(f)=\min_i\mu(f_i Q^i).\]
The fact that $Q$ is an abstract key polynomial ensures that $\mu_Q$ is a 
valuation.

Our focus in this article revolves around the following conjecture.

\begin{thcon}
Consider a valuation $\mu$ on $K[X]$ and $Q\in K[X]$ an abstract key polynomial for
$\mu$, such that $\epsilon_\mu(Q)\in\Q\otimes\mu(K^\times)$.
\begin{enumerate}
    \item There exists a subset $\Delta(Q)\subseteq\overline{K}$ that we shall call
    a \emph{diskoid}, such that for every $f\in K[X]$, the minimum over 
    $\Delta(Q)$, \ie $\min_{x\in\Delta(Q)}\mu(f(x))$ exists and is equal to the
    the truncation of $\mu$ along $Q$:

    \[\forall f\in K[X],\,\mu_Q(f)=\min_{x\in\Delta(Q)}\mu(f(x)).\]
    
    Furthermore, $\Delta(Q)$ is a finite union of balls.
    \item There is a bijective correspondence between residually transcendental
    valuations $\mu$ over $K[X]$, which can be given by a key polynomial $Q$, and
    diskoids $\Delta=\Delta(Q)$.
\end{enumerate}
\end{thcon}

What is remarkable here is that several key polynomials may induce the same
truncated valuations, but they would also induce the same diskoid $\Delta$, just 
like an open ball in a non-archimedian field may have all of its points as centres,
yet induce the same valuation. We will prove the conjecture when $K$ is henselian
or $\mu|_K$ is of rank one.

The notion of diskoids mirrors in some sense the notion of maximal divisorial set
for divisorial sets \cite{Is}.

\medskip

A different approach to the extension problem came from a series of articles by
Alexandru, Popescu and Z\u{a}h\u{a}rescu in the 90s \cite{AP}, \cite{APZ1}, 
\cite{APZ2}, \cite{APZ3}. Their idea is to extend the problem to $\overline{K}$,
the algebraic closure of $K$. The problem is simpler to study over 
$\overline{K}[X]$ since all irreducible polynomials (among which, the key
polynomials) are of degree $1$. The authors coined the concept of minimal pairs in
order to classify them. Consider a couple $(a,\delta)$ with $a\in\overline{K}$ and
$\delta=\overline{\mu}(X-a)$, where $\overline{\mu}$ stands for an arbitrary 
extension of $\mu$ to $\overline{K}[X]$. It will be called a minimal pair if
$\deg_K(a)$ is minimal among the $\deg_K(c)$, where $c\in\overline{K}$ such that
$\overline{\mu}(X-c)=\delta$. This allows to define a new valuation on
$\overline{K}[X]$ as follows
\[\overline{\mu}_{a,\delta}(f)=\min\{\overline{\mu}(a_i)+i\delta\},
\quad f=a_0+a_1(X-a)+\ldots+a_n(X-a)^n,\ a_i\in\overline{K}.\]
We then can consider the trace of $\overline{\mu}_{a,\delta}$,
by pulling it back down to $K[X]\subseteq\overline{K}[X]$. Minimal pairs have
used to study residually transcendental extensions and helped investigating
invariants associated to elements of $\overline{K}$. This work is
promoted especially by S. Khanduja (see for instance \cite{Kh}).

\medskip

In this article we show that the valuations given by minimal pairs are in fact the
same as the ones given by abstract key polynomials. Indeed we show that the trace
of the valuation given by a minimal pair over $K[X]$, is a truncated valuation of
$\mu$. We can already find a correspondence between minimal pairs and abstract key
polynomials in the work of Novacoski \cite{Nov}. In order to get closer to our
geometric interpretation we need to deepen this correspondence. We do this in our
\cref{theorem}. The correspondence can thus be stated as follows

\begin{ththm}
\begin{enumerate}
    \item Consider $Q$ the minimal polynomial of $a$, such that 
    $(a,\delta)$ is a minimal pair, for some
    $\delta\in\overline{\mu}(\overline{K}[X]^\times)$. Then
    $Q$ is an abstract key polynomial. Conversely, any abstract key
    polynomial $Q$ has an optimising root $a$, such that
    $(a,\overline{\mu}(X-a))$ is a minimal pair.
    \item Furthermore, one has an equality between valuations
    \[\overline{\mu}_{a,\delta}\big|_{K[X]}=\mu_Q.\]
\end{enumerate}
\end{ththm}

The first part of this theorem is the main goal of \cite{Nov}, the second
part is our contribution to the theory.
Similar results have been stated and proven in \cite[Theorem 5.1]{PP}.
However they are made in the context of key polynomials in the sense of
MacLane and Vaquié whereas we work in the context of abstract key
polynomials. Furthermore, they provide a proof only for augmented
valuations of residually transcendental valuations and we do not need these
assumptions. One could imagine that we could prove our statement with the
help of the work of Vaquié's extension of MacLane's theory, however we
would still need a similar result as \cite[Theorem 5.1]{PP} but for
limit-augmented valuations. Most of this work concerning these matters have
been undertaken by M. Vaquié in \cite{Va4}. Our proof will make use of the
structure of the graded algebras of a valuation which in our opinion, is an
approach worth noting. We will also give a couple of applications that we 
find noteworthy (see \cref{appli}, \cref{modif} and \cref{result}).

We also generalise \cite[Proposition 3.1]{Nov}, by using a type of Newton
polygon. This comes from our desire to encode the root configuration of a
polynomial in the Newton polygon. The first slope of this polygon will be
the $\delta$ invariant in the work surrounding minimal pairs.

\medskip

We sum up the contributions with their contributors in the following
synoptic diagram

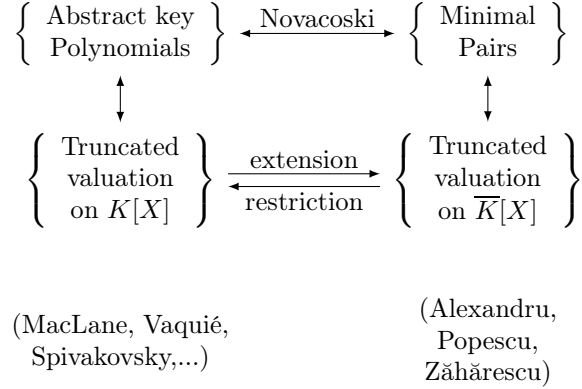
\begin{center}
    \begin{tikzcd}
    \left\{\begin{array}{c}
        \text{Abstract key}\\
        \text{Polynomials}
    \end{array}\right\}\ar[rr,leftrightarrow,"\text{\normalsize Novacoski}"]
                        \ar[d,leftrightarrow]&&
            \left\{\begin{array}{c}
                \text{Minimal}\\
                \text{Pairs}
            \end{array}\right\}\ar[d,leftrightarrow]\\
    \left\{\begin{array}{c}
        \text{Truncated}\\
        \text{valuation}\\
        \text{on }K[X]
    \end{array}\right\}\ar[rr,shift left,"\text{\normalsize extension}"]&&
            \left\{\begin{array}{c}
                \text{Truncated}\\
                \text{valuation}\\
                \text{on }\overline{K}[X]
            \end{array}\right\}\ar[ll,shift left,"\text{\normalsize restriction}"]\\
    \begin{array}{c}
        \text{(MacLane, Vaquié,}\\
        \text{Spivakovsky,...)}
    \end{array}&&
               \begin{array}{c}
                    \text{(Alexandru,}\\
                    \text{Popescu,}\\
                    \text{Z\u{a}h\u{a}rescu)}
                \end{array} 
    \end{tikzcd}
\end{center}

Our work concerns the bottom part, the extension-restriction arrows.

\medskip

Our last section is concerned with establishing the geometric
interpretation of residually transcendental valuations.
To start off, valuations given by minimal pairs can be interpreted
as minimal values of polynomials $f\in K[X]$ attained on balls. Indeed
take a minimal pair $(a,\delta)$, so that we can consider the ball
\[D(a,\delta)=\{b\in\overline{K};\:\overline{\nu}(b-a)\geqslant\delta\}.\]
This type of set is not reduced to a single point if 
$\text{rk}(\mu)=\text{rk}(\nu)$. Indeed if for instance
$\text{rk}(\mu)=\text{rk}(\nu)+1$, have $\delta$ be strictly in the last
convex subgroup of $\mu(K(X)^*)$, \ie let $\delta$ be bigger than the value
group of $K$. Then for any $x\in\overline{K}^*,\,\delta>\nu(x)$. In this
case
\[D(a,\delta)=\{a\}.\]
In order to avoid such situations we will restrict our focus to residually 
transcendental valuations. In this context
$\text{rk}(\mu)=\text{rk}(\nu)$ and thus $D(a,\delta)$ is non-reduced to a
point, if $\delta<\infty$. Then one considers a valuation 
$\overline{\nu}_{D(a,\delta)}$ defined as follows
\[\overline{\nu}_{D(a,\delta)}(f)=\min_{x\in D(a,\delta)}
\overline{\nu}(f(x)).\]
We can prove that the minimum is attained and explicit. Indeed we show that
\[\overline{\nu}_{D(a,\delta)}=\overline{\nu}_{a,\delta}.\]
Furthermore this allows for a clear bijection between residually
transcendental valuations and balls. Two minimal pairs
$(a,\delta),(a^\prime,\delta^\prime)$ may yield the same valuation, however
that is only the case if they yield the same ball
$D(a,\delta)=D(a^\prime,\delta^\prime)$. We wish to build a similar
bijection over $K$, by finding what we will call diskoids
$\Delta(Q)\subseteq\overline{K}$ for key polynomials $Q$. These will need
to verify
\begin{align*}
    \mu_Q&=\overline{\nu}_{\Delta(Q)},\text{ where }
    \overline{\nu}_{\Delta(Q)}(f)=\min_{x\in\Delta(Q)}
    \overline{\nu}(f(x))\\
    \overline{\nu}_{\Delta(Q_1)}&=\overline{\nu}_{\Delta(Q_2)}\implies
    \Delta(Q_1)=\Delta(Q_2).
\end{align*}

We can not just take $\Delta(Q)$ to be balls, since different balls can induce the
same valuation over $K$. If for instance $\mu_Q$ is the same as
$\overline{\nu}_{D(a,\delta)}\big|_{K[X]}$ and
$\sigma\in\text{Aut}_K(\overline{K})$ such that $\overline{\nu}\circ\sigma=
\overline{\nu}$, then $\mu_Q$ is also the same as the restriction to $K[X]$ of
$\overline{\nu}_{\sigma(D(a,\delta))}$. We will use bigger subsets of
$\overline{K}$.
We also want them to have a manageable shape and be easy to grasp. We will give a
reasonable candidate, called diskoids, that was studied by Julian Rüth in his PhD
thesis \cite{Ru}. We will see how diskoids decompose in simple balls and how the
absolute Galois group of $K$, \ie $G_K=\text{Gal}(K^{sep}/K)=\text{Aut}_K(\overline{K})$
acts on these balls. We wish to show that a truncated valuation (by a fixed abstract key
polynomial $Q$) is the valuation given by a diskoid (associated to $Q$) and that this
correspondence is in fact a bijection. This is what we wish to call our geometric
interpretation. We prove our statement when $(K,\nu)$ is henselian or of rank $1$.

\vspace{1cm}

\textbf{Acknowledgements.} The author wishes to thank Steven
Dale Cutkosky, Franz-Viktor Kuhlmann, Hussein Mourtada,
Bernard Teissier and the anonymous referee for the ample
discussions, corrections, remarks, suggestions and for their
interest in this work.

\section{Preliminaries}

In this paper we adopt the convention that the set $\N$ contains $0$:
$\N=\{0,1,2,3,\ldots\}$. We will write $\N^*$ for the positive integers:
$\N^*=\{1,2,3,\ldots\}$.

\subsection{Valuations and graded algebra.}

We start by defining our most basic objects, valuations. Let $R$ be a
commutative, unitary ring and $(\Gamma,+,\leqslant)$ a totally ordered
abelian group (we write $+$ for the additive law of this group and 
$\leqslant$ its order relation) to which we add an element, not in
$\Gamma$, that we denote $\infty$. We write 
$\Gamma_\infty=\Gamma\cup\{\infty\}$ and extend the addition operation $+$
and order relation $\leqslant$ so that $\infty$ plays the role of a biggest
element.

\begin{defn}
A \emph{valuation} $\nu$ on $R$, is a map $\nu\ :\ R\longrightarrow\Gamma_\infty$
satisfying
\begin{itemize}[noitemsep]
    \item[(V1)] $\forall a,b\in R,\;\nu(a\cdot b)=\nu(a)+\nu(b)$.
    \item[(V2)] $\forall a,b\in R,\;\nu(a+b)\geqslant
    \min\{\nu(a),\nu(b)\}$. This is called the ultrametric inequality.
    \item[(V3)] $\nu(1)=0$ and $\nu(0)=\infty$.
    \item[(V4)] $\nu^{-1}(\infty)=(0)$
\end{itemize}
We often write $(R,\nu)$ for a \emph{valuative pair}, \ie a ring equipped
with a valuation.
\end{defn}

\begin{rke}
\begin{enumerate}
    \item When we talk of embedding a pair $(R,\nu)$ into $(S,\mu)$, we
    mean we set an injective morphism $R\xhookrightarrow{\iota} S$, such
    that $\mu\circ\iota=\nu$.
    \item (V4) implies that $R$ is a domain. One could then extend $\nu$ to
    a valuation over $Q(R)$, the quotient field of $R$, by setting
    $\nu(a/b)=\nu(a)-\nu(b)$. It is easy to see it is well-defined.
\end{enumerate}
\end{rke}

Our study of valuations requires a powerful tool, the graded algebra of a
valuation. For $\gamma\in\nu(R\setminus\{0\})$ and define the following
groups
\begin{align*}
    \mathcal{P}_\gamma=\mathcal{P}_\gamma(R,\nu) &=
                \{a\in R\ |\ \nu(a)\geqslant\gamma\}\\
    \mathcal{P}^+_\gamma=\mathcal{P}^+_\gamma(R,\nu) &=
                \{a\in R\ |\ \nu(a)>\gamma\}.
\end{align*}

The graded ring $\text{gr}_\nu(R)$ is
\[\text{gr}_\nu(R)=\bigoplus_{\gamma\in\nu(R\setminus\{0\})}
\frac{\mathcal{P_\gamma}}{\mathcal{P}^+_\gamma}.\]
This comes equipped with a map, called the initial form, $\init_\nu\ :\ 
R\setminus\{0\}\longrightarrow\text{gr}_\nu(R)$ assigning to $a$ its
class modulo $\mathcal{P}_\gamma^+(R)$ with $\gamma=\nu(a)$. By definition,
any homogeneous element is thus represented and any initial form of any
non-zero element is non-zero. We can assign to $0$ the value $0$ in the
graded ring.

\begin{rke}
For general filtered modules (or rings, algebras etc.) there is also a
notion of initial form which fails to be a morphism in general. It may not
even be multiplicative, however when considering the graded algebra
associated to a valuation, we can still state some simple rules of
computation:
\begin{enumerate}[noitemsep]
    \item $\init_\nu(a\cdot b)=\init_\nu(a)\cdot\init_\nu(b),\;
    \forall a,b\in R$. Thus the graded ring is an integral domain.
    \item if $\nu(a)>\nu(b)$ then $\init_\nu(a+b)=\init_\nu(b)$.
    \item if $\nu(a)=\nu(b)<\nu(a+b)$ then 
        $\init_\nu(a+b)\neq\init_\nu(a)+\init_\nu(b)=0$.
    \item if $\nu(a)=\nu(b)=\nu(a+b)$ then
        $\init_\nu(a+b)=\init_\nu(a)+\init_\nu(b)$.
\end{enumerate}
1. is simply a consequence of (V1) and the rest amount to using (V2).
\end{rke}

\subsection{Numerical Invariants.}

There are many ways in which one can measure the complexity of a valuation.
We start with assigning invariants to its value group.

\begin{defn}
Consider a valuative pair $(R,\nu)$ and $K$ the fraction field of $R$.
The \emph{value group of $\nu$} written $\Phi(\nu)$ is
\[\Phi(\nu):=\nu(K^\times).\]

We define the \emph{rank of $\nu$}, written $\text{rk}(\nu)$
\[\text{rk}(\nu):=\text{ord}\{\text{convex subgroups of }\Phi(\nu)\}\]
where ord denotes the ordinal type of the set. We restrict ourselves to
the cases where indeed, the set of convex subgroups (also called
\emph{isolated subgroups} in the literature) do form a well-ordered set.
Next we assign to $\nu$ its \emph{rational rank}, written
$\text{rat.rk}(\nu)$
\[\text{rat.rk}(\nu):=\dim_\Q \Phi(\nu)\otimes_\Z \Q.\]
\end{defn}

\begin{rke}
It is quite possible that $\nu(R\setminus\{0\})$ is a lot smaller than 
$\Phi(\nu)$ (or may not even be a group), however it is clear that
$\nu(R\setminus\{0\})$ generates $\Phi(\nu)$ as a subgroup.
\end{rke}

\begin{defn}
The \emph{valuation ring} of the pair $(R,\nu)$ is a local subring of $K$,
written $\Or(\nu)$ with maximal ideal $\m(\nu)$ and residual field
$\kappa(\nu)$. We define them as follows:
\begin{align*}
    \Or(\nu) &= \{a\in K\ |\ \nu(a)\geqslant0\}\\
    \m(\nu) &= \{a\in K\ |\ \nu(a)>0\}\\
    \kappa(\nu) &= \Or(\nu)/\m(\nu).
\end{align*}
\end{defn}

By \emph{extension of valued rings}, we mean an extension of rings
$R\subseteq S$, each equipped with a respective valuation $\nu$ and $\mu$,
such that $\mu$ restricts to $\nu$ on $R$,
\[\mu\big|_{R}=\nu.\]
This is equivalent to saying that the natural inclusion map
$R\hookrightarrow S$ defines an embedding of valued pairs
$(R,\nu)\hookrightarrow(S,\mu)$. Given such an extension one can compare
value groups and residual fields. Indeed, one has
$\nu(R\setminus\{0\})\subset\mu(S\setminus\{0\})$, thus $\Phi(\nu)\subset
\Phi(\mu)$. Secondly we have $\Or(\nu)\subseteq\Or(\mu)$ and
$\Or(\nu)\cap\m(\mu)=\m(\nu)$ thus giving us a natural inclusion of
residual fields $\kappa(\nu)\subseteq\kappa(\mu)$. We will call them
\emph{value group extension} and \emph{residual extension} respectively.

\begin{prop}
We have the following inequalities given for any valuation $\nu$ on $K$ and
any valued extension $(L,\mu)$:
\begin{enumerate}
    \item $\text{rk}(\nu)\leqslant\text{rat.rk}(\nu)$.
    \item $\dim_\Q (\Phi(\mu)/\Phi(\nu))\otimes_\Z \Q+
    \text{tr.deg}_{\kappa(\nu)}\kappa(\mu)\leqslant\text{tr.deg}_K L$
\end{enumerate}
when these quantities are well-defined and finite.
The second one is called the Zariski-Abhyankar inequality (see
\cite[Ch. 6, §10, no. 3, Cor. 1]{Bo} for a proof). One could even
generalise it (see \cite[Appendix 2, Prop.~2]{ZS2} or
\cite[Theorem 1]{Ab2}).
\end{prop}

This last result allows us to say that the rank of a valuation $\mu$ on
$K[X]$, extending a valuation $\nu$ on $K$ can not jump more than once
\[\text{rat.rk}(\nu)\leqslant\text{rat.rk}(\mu)\leqslant
\text{rat.rk}(\nu)+1.\]
We will say that our (simple transcendental) extension is \emph{valuation
algebraic} if the quotient group $\Phi(\mu)/\Phi(\nu)$ is torsion and if
the residual extension $\kappa(\mu)/\kappa(\nu)$ is algebraic. Otherwise it
is \emph{valuation transcendental}. It will be called \emph{value
transcendental} if $\Phi(\mu)/\Phi(\nu)$ is not torsion and
\emph{residually transcendental} if $\kappa(\mu)/\kappa(\nu)$ is not
algebraic. The Abhyankar inequality tells us that either cases are
possible, but not both at the same time.

\section{Truncation and Abstract Key Polynomials}\label{secabkp}

Given two polynomials $f, Q\in K[X],\:\deg Q>0$, one can define by
successive euclidean divisions, the $Q$-expansion of $f$
\[f=f_0+f_1 Q+\ldots+f_n Q^n,\;\forall i,\,\deg f_i<\deg Q.\]
This expansion is unique. When given a valuation $\mu$ on $K[X]$ one can
then define the \emph{truncated map of $\mu$, with respect to $Q$}
\[\mu_Q(f)=\min_i \mu(f_iQ^i).\]

For $f\in K[X]$ and $Q,\mu$ as above, we will write
\[S_{Q,\mu}(f)=\{i\in\N\ |\ \mu(f_i Q^i)=\mu_Q(f)\}\]
or just $S_Q(f)$ for short, when the $\mu$ will be fixed once and for all.
We also define $d_{Q,\mu}(f)=d_Q(f):=\max S_{Q,\mu}(f)$.

The truncated map is not always a valuation. It still is a map extending
the valuation $\nu=\mu\big|_K$ and it still verifies (V2) and (V3) (we
assume that our valuation has trivial support and thus, $\mu_Q$ also checks
(V4)). It may however fail to verify (V1). See \cite[Example 2.5]{NS} for a
counter-example. There is a natural condition for which it is a valuation,
that is if $Q$ is an abstract key polynomial. In that case we will talk
about \emph{truncated valuation in} $Q$. We will need to define the
$\epsilon$-level or $\epsilon$ factor of a polynomial $f$.

\begin{defn}
For any valuation $\mu$ on $K[X]$ and polynomial $f\in K[X]$, just as
above, we define
\[\epsilon_\mu(f)=\max_{i\geqslant 1}\frac{\mu(f)-\mu(\partial_i f)}{i}.\]
We will call it the $\epsilon$ factor of $f$.
For our fixed $\mu$ we will simply write $\epsilon(f)=\epsilon_\mu(f)$ and
when we will truncate by $Q$, we will write
$\epsilon_Q(f)=\epsilon_{\mu_Q}(f)$.
\end{defn}

Here $\partial_i$ represents the formal Hasse-Schmidt derivative on $K[X]$.
This is an operator that can be defined by means of the Taylor expansion of
polynomials in two variables
\[f(X+Y)=\sum_{i\geqslant 0} \partial_i f(X) Y^i.\]

By multiplying together the Taylor expansions of two polynomials $f,g\in
K[X]$, we see that these Hasse-Schmidt derivatives satisfy the Leibniz rule
\[\partial_i(f\cdot g)=\sum_{i=j+k}\partial_j f\cdot\partial_k g.\]

Furthermore we can compose Hasse-Schmidt derivatives. By expanding
$f(X+Y+Z)$ in two different ways, we can show that
\[\partial_i\circ\partial_j=\binom{i+j}{j}\partial_{i+j}.\]

We now define abstract key polynomials.

\begin{defn}\cite[Definition 11]{DMS}
Let $Q\in K[X]$ be a monic polynomial. We say that $Q$ is an
\emph{abstract key polynomial for $\mu$} (we will abbreviate ABKP) if,
for any $f\in K[X]$
\[\deg f<\deg Q\implies\epsilon_\mu(f)<\epsilon_\mu(Q).\]
\end{defn}

As basic examples, any degree one polynomials are abstract key polynomials
according to this definition. Less obvious examples are the key polynomials
given by MacLane-Vaquié's key polynomials. See \cite[Section 3]{DMS}.

Now let us go through some basic properties the ABKPs verify.

\begin{prop}\cite[Prop. 2.4 (ii)]{NS}
Abstract key polynomials are irreducible.
\end{prop}

\begin{prop}\cite[Proposition 13]{DMS}
Let $P_1,\ldots,P_t\in K[X]$ be polynomials of degree $<\deg Q$
(assume $t\geqslant 2$). If we set the following euclidean division
$\prod_{i=1}^t P_i=qQ+r$ in $K[X]$, with $\deg r<\deg Q$, then
\[\mu\left(\prod_{i=1}^t P_i\right)=\mu(r)<\mu_Q(qQ)\leqslant\mu(qQ).\]
\end{prop}

Remark that in the above proposition, $r$ can not be $0$, since $Q$ is 
irreducible, thus $\prod_{i=1}^t P_i$ is prime with $Q$, since each $P_i$
is prime with $Q$.

\begin{corr}\cite[Proposition 15]{DMS}
If $Q$ is an ABKP, then $\mu_Q$ is a valuation.
\end{corr}

This allows us to define the graded algebra of $\mu_Q$. We wish to
understand its structure, but we first need the following.

\begin{corr}\cite[Remark 16]{DMS}
Let $\alpha=\deg Q$ and define
\[G_{<\alpha}=\sum_{\deg f<\deg Q}\text{gr}_\nu(K)\cdot\init_{\mu_Q}(f)
\subseteq\text{gr}_{\mu_Q}(K[X])\]
then this $\text{gr}_\nu(K)$-module is stable under multiplication and thus
is an algebra.
\end{corr}

\begin{rke}\label{sum}
We will see later (see \cref{modif}) that the condition $\deg f<\deg Q$
under the sum can be replaced by $\epsilon_\mu(f)<\epsilon(Q)$.
\end{rke}

This subring $G_{<\alpha}$ will play the role of coefficients in the
following proposition.

\begin{prop}\cite[Remark 16]{DMS}
The graded ring of $\mu_Q$ on $K[X]$ has a simple polynomial structure.
More precisely
\[\text{gr}_{\mu_Q}(K[X])=G_{<\alpha}\left[\init_{\mu_Q}Q\right]\]
with $\init_{\mu_Q}Q$ transcendental over $G_{<\alpha}$.
\end{prop}

Let us just describe the action of the initial form here. Consider a
polynomial $f\in K[X]$ and its $Q$-expansion
\[f=f_0+f_1 Q+\ldots+f_n Q^n,\;\forall i,\,\deg f_i<\deg Q\]
so that 
\[\init_{\mu_Q}(f)=\sum_{j\in S_Q(f)}
\init_{\mu_Q}(f)\init_{\mu_Q}(Q)^j\]
thanks to the rules of computation with initial forms. Thus
$\init_{\mu_Q}(f)\in G_{<\alpha}[\init_{\mu_Q}(Q)]$ and
\[\deg_{\init_{\mu_Q}}\left(\init_{\mu_Q}(f)\right)=
\max S_{Q,\mu}(f).\]

Truncation now gives us a valuation. Furthermore, we have that
a truncation will bound the $\epsilon$-level associated to the
truncation. We state the main result concerning this numerical
properties, but we will postpone its proof (see \cref{result})
as we would have by then introduced the tools to
prove it in a fairly straightforward way.

\begin{prop}\cite[Lemma 17]{DMS}\label{num}
For any polynomial $f\in K[X]$
\[\epsilon_{\mu_Q}(f)\leqslant\epsilon_{\mu_Q}(Q)=
\epsilon_\mu(Q).\]

Furthermore we have a necessary and sufficient condition for 
equality:
\[\epsilon_{\mu_Q}(f)=\epsilon_{\mu_Q}(Q)\iff
S_{Q,\mu}(f)\neq\{0\}.\]
\end{prop}

\section{Newton polygons}

In this section we will present a result concerning the values 
$\overline{\mu}(X-\lambda)$ where $\lambda$ runs through the
roots of a polynomial $f$. We will call this data the root
configuration of $f$. We will use a certain type of Newton
polygons. The first slope of the polygon will be the $\epsilon$
factor. It is our hope that the whole slope data will lead us
to a better understanding of the diskoids decomposition and the
action of the absolute Galois group of $K$ on them.

\subsection{Defining our Newton polygon.}
Classical Newton polygons take a field $F$ equipped an
ultrametric absolute value $|.|$ (or a valuation $V$ of rank
$1$, so that one can suppose that $V(F^\times)\subseteq\R$). If
however we want to work with fields of higher ranks, one needs
to work out what convex sets are in $\R\times\Phi$ for any
ordered abelian group $\Phi$.
Our presentation here takes many components from Vaquié's
presentation in \cite{Va2}. We will consider the base changes
$\Phi_\Q=\Phi\otimes_\Z\Q$ and $\Phi_\R=\Phi\otimes_\Z\R$.
Since $\Phi$ has no torsion, the canonical maps
\begin{center}
    \begin{tikzcd}
    \Phi\ar[r]&\Phi_\Q\ar[d,hook]&\gamma\ar[r,mapsto]&
        \gamma\otimes 1\\
    \Phi\ar[r]&\Phi_\R&\gamma\ar[r,mapsto]&\gamma\otimes 1\\
    \end{tikzcd}
\end{center}
are all injective, thus we can consider $\Phi$ as a subgroup of
$\Phi_\Q$ or $\Phi_\R$.

We define a \emph{line} to be a subset 
$L\subseteq\R\times\Phi_\R$ defined by a linear or affine
equation
\[L=L_{q,\alpha,\beta}=\{(x,\gamma)\in\R\times\Phi_\R\ |\
q\gamma+\alpha x+\beta=0\}\]
for some fixed values $q\in\R,\;\alpha,\beta\in\Phi_\R$. The
slope of $L$ is given by $s(L)=s(L_{q,\alpha,\beta})=\alpha/q$
whenever $q\neq0$. This definition of slope is not classical
and corresponds to the negative of the natural slope of a line.

There will always be a single line passing through two fixed
and distinct points $P_1,P_2$, that we will denote $(P_1,P_2)$.
Any line defines two \emph{half-spaces}
\begin{align*}
    H^L_\geq &= \{(x,\gamma)\in\R\times\Phi_\R\
                    |\ q\gamma+\alpha x+\beta\geq 0\}\\
    H^L_\leq &= \{(x,\gamma)\in\R\times\Phi_\R\
                    |\ q\gamma+\alpha x+\beta\leq 0\}.\\
\end{align*}
For any subset $A\subseteq\R\times\Phi_\R$, we define its
\emph{convex hull} by
\[\text{Conv}(A)=\bigcap_{\substack{H\text{ half-space } \\
A\subseteq H}}H\]
\ie the intersection of half-spaces containing $A$. A set $A$
is considered to be convex if $\text{Conv}(A)=A$.
For any set $A$, we define its \emph{faces} to be subsets
$F=\text{Conv}(A)\cap L$ where $L$ is a line in 
$\R\times\Phi_\R$, satisfying
\begin{itemize}[noitemsep]
    \item[-] $\text{Conv}(A)$ is contained in one of the
        half-spaces $H^L_\geq$ or $H^L_\leq$.
    \item[-] $F=\text{Conv}(A)\cap L$ contains at least two
        points.
\end{itemize}
We will say that $L$ \emph{supports the face $F$}. The slope
$s(F)$ of a face $F$ of $A$ will simply be $s(L)$ where $L$
supports $F$. Usually a Newton polygon is constructed for
finite sets $X=\{(k,\gamma_k),\:0\leqslant k\leqslant m\}$. We
will write its \emph{Newton polygon} as
\[PN(X)=\text{Conv}(\{(x,\delta)\ |\ \exists(x,\gamma)\in X,\;
\delta\geq\gamma\})=\text{Conv}\left((\{0\}\times
\Phi_{\geq 0})+X\right)\]
where $\Phi_{\geq0}=\{\gamma\in\Phi_\R,\;\gamma\geq0\}$ and 
$A+B$ is the Minkowski sum of two subsets of $\R\times\Phi_\R$.
The bottom boundary of $PN(X)$ is then a finite polygonal line,
thus describing it will be equivalent to giving the following

\begin{itemize}[noitemsep]
    \item[-] a finite sequence of non-negative integers $0=a_0<a_1<\ldots<a_r=m$
    (the abscissa or $x$-coordinates of the vertices of the polygonal line),
    \item[-] a finite set of values in $\Phi_\R:\;\epsilon_1>\ldots>\epsilon_r$
    (the slopes of the segments forming the polygonal line)
\end{itemize}
verifying

\begin{enumerate}[noitemsep]
    \item $\forall k,t,\, 0\leqslant k\leqslant m,1\leqslant 
        t\leqslant r$ one has
    \[\gamma_k+k\epsilon_t\geqslant\gamma_{a_{t-1}}+a_{t-1}
        \epsilon_t=\gamma_{a_t}+a_t\epsilon_t.\]
    \item if $k<a_{t-1}$ or $k>a_t$, then the above inequality 
    is strict (this should be interpreted as $(a_t,\gamma_t)$
    being the points where the polygonal line turns).
\end{enumerate}

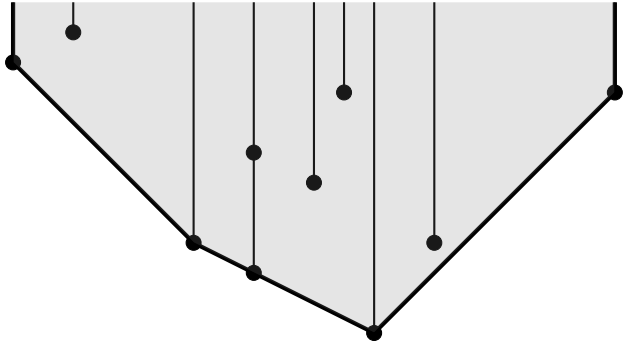
\begin{figure}[ht]
    \centering
    \begin{tikzpicture}[scale=0.4]
    \node[draw,circle,inner sep=2pt,fill] at (0,9) {};
    \node[draw,circle,inner sep=2pt,fill] at (2,10) {};
    \node[draw,circle,inner sep=2pt,fill] at (6,3) {};
    \node[draw,circle,inner sep=2pt,fill] at (8,2) {};
    \node[draw,circle,inner sep=2pt,fill] at (8,6) {};
    \node[draw,circle,inner sep=2pt,fill] at (10,5) {};
    \node[draw,circle,inner sep=2pt,fill] at (11,8) {};
    \node[draw,circle,inner sep=2pt,fill] at (12,0) {};
    \node[draw,circle,inner sep=2pt,fill] at (14,3) {};
    \node[draw,circle,inner sep=2pt,fill] at (20,8) {};
    
    \draw[ultra thick, black] (0,11) -- (0,9);
    \draw[ultra thick, black] (0,9) -- (6,3);
    \draw[ultra thick, black] (6,3) -- (12,0);
    \draw[ultra thick, black] (12,0) -- (20,8);
    \draw[ultra thick, black] (20,8) -- (20,11);
    
    \draw[thick, black] (2,10) -- (2,11);
    \draw[thick, black] (6,3) -- (6,11);
    \draw[thick, black] (8,2) -- (8,11);
    \draw[thick, black] (10,5) -- (10,11);
    \draw[thick, black] (11,8) -- (11,11);
    \draw[thick, black] (12,0) -- (12,11);
    \draw[thick, black] (14,3) -- (14,11);
    
    \filldraw[fill=gray, fill opacity=0.2, draw=black] (0,11) -- (0,9) -- (6,3) -- (12,0) -- (20,8) -- (20,11);
    \end{tikzpicture}
    \caption{Example of a Newton polygon of a finite subset of $\R\times\R$.}
\end{figure}

The points $A_t=(a_t,\gamma_{a_t})$ are the \emph{vertices} of
the polygon and the \emph{faces} are simply the segments
$[A_{t-1},A_t],\,1\leqslant t\leqslant r$ that have slope
$\epsilon_t$. We will also call $a_t-a_{t-1}$ the \emph{length}
of this face.

\medskip

Let us now fix a valuation $\mu$ over $K[X]$ with values in
$\Phi$. For any $f\in K[X]$ define
\[PN(f,\mu)=PN(X(f,\mu))\text{ where }
X(f,\mu)=\{(i,\mu(\partial_i f))\ |\ i=0,\ldots,\deg(f)\}.\]

\begin{defn}
For any valuation $\mu$ on $K[X]$ and polynomial $f\in K[X]$,
we define the \emph{slope data} $f$, as the sequence of
vertices of $PN(f,\mu)$
\[\big((a_1,\epsilon_1),\ldots,(a_r,\epsilon_r)\big)\]
as defined above.
\end{defn}

In the next subsection, we relate the finite data given by our
version of the Newton polygon of $f$ to some information given
by the configuration of roots of $f$.

\subsection{Root configurations.}

In order to define what we mean by configuration of roots, we 
need to extend our initial valuation $\mu$ to $\overline{K}[X]$
where $\overline{K}$ is the algebraic closure of $K$. Decompose
then $f$ into linear factors
\[f=\alpha\prod_{i=1}^m(X-\lambda_i).\]
What will be of interest to us will be the values 
$\overline{\mu}(X-\lambda_i)$ and how many times one such value
is repeated. We can fix an indexation of the $a_i$ so that we
have the following

\begin{align*}
    \overline{\mu}(X-\lambda_1)=\ldots=
        \overline{\mu}(X-\lambda_{l_1}) &
        >\overline{\mu}(X-\lambda_{l_1+1})=\ldots=
            \overline{\mu}(X-\lambda_{l_1+l_2})\\
        &\vdots\\
        &>\overline{\mu}(X-\lambda_{l_1+\ldots+l_{s-1}+1})=\ldots=
            \overline{\mu}(X-\lambda_{l_1+\ldots+l_s}).
\end{align*}

We will write $\overline{\mu}(X-\lambda_{l_1+\ldots+l_t})=
\delta_t$.

\begin{defn}
Set $\overline{\mu}$ and $f$ as above. Adopting the previous 
notations, we define the \emph{root configuration of} $f$ as
the finite sequence
\[\big((l_1,\delta_1),\ldots,(l_s,\delta_s)\big).\]
\end{defn}

Our main result in this section says that for any polynomial
$f$, its root configuration and slope data are equivalent.

\begin{thm}\label{thm}
Fix $\overline{\mu}$ a valuation over $\overline{K}[X]$ and
$f\in K[X]$. Then the root configuration of
$f,\;\big((l_t,\delta_t),t=1,\ldots,s\big)$ is encoded in its
slope data $\big((a_t,\epsilon_t),t=1,\ldots,r\big)$ as
follows:
\begin{align*}
    s &= r\\
    a_t-a_{t-1} &= l_t,\: t=1,\ldots,r\\
    \epsilon_t &= \delta_t,\: t=1,\ldots,r
\end{align*}
where we have fixed $a_0=0$.
\end{thm}

\begin{rke}
One can be surprised by the fact that even though the 
definition of the $\delta_t$ needs us to choose an extension
$\overline{\mu}$ of $\mu$ to $\overline{K}$, their value is
ultimately independent of this choice.
\end{rke}

To prove \Cref{thm} one can use a classic lemma concerning
Newton polygons:

\begin{lm}\label{Kob}
Let $(F,V)$ be any valued field and $p(T)\in F[T]$ any 
polynomial, whose roots are in $F$. We write
\[p(T)=c\prod_{k=1}^N(T-c_k)=\sum_{l=0}^N b_l T^l
\quad c,c_k,b_l\in F\]
and consider the points $\{(l,V(b_l)),l=0,\ldots,N\}$. If
$\zeta$ is a slope of its Newton polygon of length $\ell$, it
follows that precisely $\ell$ of the roots $c_k$ have valuation
$\zeta$.
\end{lm}

A proof of \Cref{Kob} can be found in \cite[Ch. IV, § 3]{Kbtz}
in the rank $1$ case, but the proof works verbatim for $V$ of
arbitrary rank as well. We can apply it to prove \Cref{thm}, by
setting $(F,V)=(\overline{K}(X),\overline{\mu})$ and
\[p(T)=f(X+T)=
\sum_{i=0}^m \partial_if(X) T^i\in \overline{K}(X)[T].\]
Indeed, the roots of $p(T)$, when considered as a polynomial of
coefficients in $\overline{K}(X)$ are $\lambda_i-X,
i=1,\ldots,m$ \footnote{This short proof was suggested by
Franz-Viktor Kuhlmann in a private communication.}.

\subsection{The $\mathbf{\delta}$ invariant.}
In \cite{Nov}, another quantity is defined in parallel to the
$\epsilon$ factor.

\begin{defn}
For any valuation $\overline{\mu}$ of $\overline{K}[X]$ and any
polynomial $f\in\overline{K}[X]$, we define $\delta$ as follows
\[\delta_{\overline{\mu}}(f):=\max\left\{\overline{\mu}(X-c);\:
c\text{ root of }f\right\}.\]
We will often abbreviate to $\delta(f)$ when the valuation 
$\overline{\mu}$ is implicit.
\end{defn}

Considering the root configuration
$\big((l_t,\delta_t),t=1,\ldots,r\big)$ of $f$, we obviously
have by definition
\[\delta(f)=\delta_1.\]

On the other hand, if $\big((l_t,\delta_t),t=1,\ldots,r\big)$
is the slope data of $f$, then we can write down
\[\epsilon(f)=\epsilon_1.\]
Thus \cite[Prop. 3.1]{Nov} is a consequence of our \cref{thm}
for $t=1$.

\begin{rke}
We should observe that even though the choice $\overline{\mu}$
is arbitrary, the $\delta_{\overline{\mu}}(f)$, when specifying
$f\in K[X]$, will only depend on $\mu$.
\end{rke}

\section{Minimal Pairs}

In this section, we wish to build valuations over $K[X]$ by bringing one 
down from $\overline{K}[X]$. We start with $(K(X),\mu)$ and extend the
valuation to $(\overline{K}(X),\overline{\mu})$. This can be represented in
a simple diagram of extensions of valuative pairs

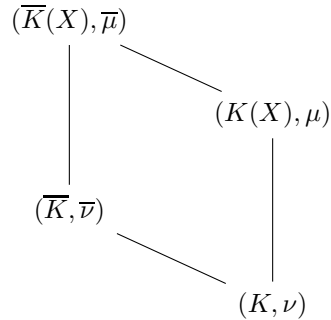
\begin{center}
    \begin{tikzcd}
    (\overline{K}(X),\overline{\mu})&\\
    &(K(X),\mu)\ar[ul,no head]\\
    (\overline{K},\overline{\nu})\ar[uu,no head]&\\
    &(K,\nu)\ar[ul,no head]\ar[uu,no head]
    \end{tikzcd}
\end{center}

\begin{rke}
We can show that $\text{Aut}_K(L)\cong\text{Aut}_{K(X)}(L(X))$ for any 
algebraic field extension $L/K$. Indeed the isomorphism is given by simply
taking a $K(X)$-automorphism $\sigma$ of $L(X)$ and restricting it to $L$.
In order to see that this mapping is well-defined, it is enough to show that 
if $\varphi\in\text{Aut}_{K(X)}(L(X))$ and $\alpha\in L$ then
$\varphi(\alpha)\in L$. We know that
$\alpha$ is algebraic over $K$, thus $\varphi(\alpha)$ is
algebraic as well: there is a $P\in K[X]$ such that 
$P(\alpha)=0$. If $\varphi(\alpha)\in L(X)\setminus L$, then
$\varphi(\alpha)$ is transcendental over $L$ and thus,
transcendental over $K$, which is a contradiction.

To see how to extend an element $\sigma\in\text{Aut}_K(L)$ to an element
$\Tilde{\sigma}\in\text{Aut}_{K(X)}(L(X))$, take any $f=\sum_k a_k X^k\in 
L[X]$ and set

\[\Tilde{\sigma}(f)=\sum_k \sigma(a_k)X^k\]
and this can then be extended to $L(X)$ by setting $\Tilde{\sigma}(f/g)=
\Tilde{\sigma}(a)/\Tilde{\sigma}(b)$. This defines the reverse mapping to the
one defined above.

From now on we will identify elements from the two groups.
\end{rke}

\begin{lm}\label{lm1}
For any situation where we extend $\nu$ to $\overline{\nu}$ over $\overline{K}$ and
$\nu$ to $\mu$ over $K[X]$, there is a common extension $\overline{\mu}$, of both
$\mu$ and $\overline{\nu}$.
\end{lm}

\begin{proof}
We begin by taking any extension $\mu^\prime$ of $\mu$ to
$\overline{K}(X)$. Since $\overline{K}/K$ is normal,
$\text{Aut}_K(\overline{K})$ acts transitively on the
valuations of $\overline{K}$ extending $\nu$, so that there is
a $\sigma\in\text{Aut}_{K(X)}(\overline{K}(X))$ such that
$\mu^\prime\circ\sigma\big|_{\overline{K}}=\overline{\nu}$.
\end{proof}

If we wish to construct a valuation on $\overline{K}[X]$, we
can consider $a\in\overline{K}$, $\delta\in\Phi$, where $\Phi$
is a group containing $\Phi(\overline{\nu})$, and define
\[\overline{\nu}_{a,\delta}\left(\sum_i a_i(X-a)^i\right)=
\min_i\{\overline{\nu}(a_i)+i\delta\}.\]
We say that $(a,\delta)$ is a \emph{pair of definition} for
a valuation $\overline{\mu}$ if $\overline{\mu}=
\overline{\nu}_{a,\delta}$.

\begin{rke}
\begin{enumerate}[noitemsep]
    \item One observes that $a_i=\partial_if(a)$.
    \item If $\delta=\overline{\mu}(X-a)$ then
    $\overline{\nu}_{a,\delta}$ is a truncation:
    $\overline{\nu}_{a,\delta}=\overline{\mu}_{X-a}$.
\end{enumerate}
\end{rke}

Several pairs of definition can yield the same valuation. We
characterise these situations in the lemma below.

\begin{lm}\label{pairs}
Two pairs, $(a,\delta)$ and $(a',\delta')$ define the same
valuation if and only if the following two conditions hold
\begin{enumerate}[noitemsep]
    \item $\delta=\delta'$
    \item $\overline{\nu}(a-a')\geqslant\delta$
\end{enumerate}
\end{lm}

There is a proof of this fact in the case of residually
transcendental extensions in \cite{AP}, but we give a proof
in the general case. 

\begin{proof}
Consider two pairs $(a,\delta),(a',\delta')$ that define the 
same valuation. Then we have by definition
\[\delta'=\overline{\nu}_{a',\delta'}(X-a')=
\overline{\nu}_{a,\delta}(X-a')=
\min\{\delta,\overline{\nu}(a-a')\}\]
and by a symmetric argument we obtain
\[\delta=\min\{\delta',\overline{\nu}(a-a')\}\]
thus $\delta'=\delta$ and 
$\overline{\nu}(a-a')\geqslant\delta$.

Conversely, consider pairs that verify the two conditions of
the lemma and let us show that the valuations they define are
equal. It is clear that they agree on $\overline{K}$ and it
is sufficient to show they agree on polynomials of type
$X-b,\, b\in\overline{K}$ to show that they are equal.
Indeed, each polynomials factors as products of such simple 
degree 1 polynomials, since $\overline{K}$ is algebraically
closed. We have
\begin{align*}
    \overline{\nu}_{a,\delta}(X-b)&=\min\{\delta,
        \overline{\nu}(a-b)\}\\
    \overline{\nu}_{a',\delta}(X-b)&=\min\{\delta',
        \overline{\nu}(a'-b)\}.
\end{align*}
If $\overline{\nu}(a-b)\geqslant\delta$ then
$\overline{\nu}_{a,\delta}(X-b)=\delta$, but
$\overline{\nu}(a'-b)=\overline{\nu}(a'-a+a-b)\geqslant
\min\{\overline{\nu}(a'-a),\overline{\nu}(a-b)\}\geqslant
\delta$ so that
\[\overline{\nu}_{a',\delta}(X-b)=\delta=
\overline{\nu}_{a,\delta}(X-b).\]
If $\overline{\nu}(a-b)<\delta$ then
$\overline{\nu}(a'-b)=\overline{\nu}(a'-a+a-b)=
\overline{\nu}(a-b)$ thus
\[\overline{\nu}_{a,\delta}(X-b)=\overline{\nu}(a-b)=
\overline{\nu}(a'-b)=\overline{\nu}_{a',\delta}(X-b).\]
\end{proof}

One can choose, among the $\{a'\in\overline{K};\:
\overline{\nu}(a'-a)\geqslant\delta\}$ one of minimal degree.
This leads us to the natural definition of minimal pairs of
definition.

\begin{defn}
We will say that $(a,\delta)$ is a \emph{minimal pair of
definition for a valuation} $\overline{\mu}$, if
$\overline{\mu}=\overline{\nu}_{a,\delta})$, so that
$(a,\delta$ is a pair of definition for $\overline{\mu}$ and
$\deg_K(a)$ is minimal among the pairs that define it.
\end{defn}

If we are given an arbitrary $\overline{\mu}$ one can try to
approximate it by setting $\delta=\overline{\mu}(X-a)$ so
that $\overline{\nu}_{(a,\delta)}\leqslant\overline{\mu}$.
The definition of minimal pair consists of a couple formed by
an element $a$ and the value $\delta=\overline{\mu}(X-a)$, so 
that $\deg_K(a)$ is minimal among the centres defining
$\overline{\nu}_{a,\delta}$.

\begin{defn}
We say that $(a,\delta)$ is a \emph{minimal pair} for 
$\overline{\mu}$ when
\begin{enumerate}[noitemsep]
    \item $\overline{\mu}(X-a)=\delta$.
    \item for any $b\in\overline{K},\;\deg_K(b)<\deg(a)
    \implies\overline{\nu}(a-b)<\delta$.
\end{enumerate}
\end{defn}

\begin{rke}
\begin{enumerate}[noitemsep]
    \item Condition 2 of Definition 4.6 above is also 
    equivalent to the contrapositive:
    \[2'.\:\forall b\in\overline{K},\;\overline{\nu}(a-b)
    \geqslant\delta\implies\deg_K(b)\geqslant\deg(a).\]
    \item A minimal pair of definition characterises the 
    valuation we are studying and allows for a direct way of
    computing it. A minimal pair does not allow to compute the
    values of a polynomial: a valuation $\overline{\mu}$ may
    have a minimal pair $(a,\delta)$, however it may not be a
    minimal pair of definition for $\overline{\mu}$ as
    $\overline{\nu}_{a,\delta}$ may be different from
    $\overline{\mu}$.
    \item Minimal pairs (of definition) were defined in order 
    to study residually transcendental valuation extensions,
    abbreviated RT extension (\ie the extension
    $\kappa(\mu)/\kappa(\nu)$ is of transcendence degree $1$).
\end{enumerate}
\end{rke}

RT extensions are characterised in the theorem below.

\begin{thm}\cite[Proposition 1, 2 and 3]{AP}\label{rt}
The following are equivalent
\begin{enumerate}[noitemsep]
    \item $\mu/\nu$ is an RT extension.
    \item $\overline{\mu}/\overline{\nu}$ is an RT extension.
    \item $\exists a\in\overline{K},\delta\in\Phi(\overline{\nu})$ such that
    $\overline{\mu}=\overline{\nu}_{a,\delta}$.
    \item $\Phi(\overline{\mu})=\Phi(\overline{\nu})$ and
    $\overline{\mu}(X-\overline{K}):=\{\overline{\mu}(X-b),\ b\in\overline{K}\}$
    has a maximal element.
\end{enumerate}
Furthermore, if one of these conditions is verified, we then have:
\[\max \overline{\mu}(X-\overline{K})=\delta.\]
\end{thm}

\begin{rke}
Item 3 of \cref{rt} will be refined for the extension
$(K[X],\mu)/(K,\nu)$ once we have established \cref{theorem}
(see \cref{appli}).
\end{rke}

When given a minimal pair $(a,\delta)$, one can compute
$\mu(f)$ for polynomials such that $\deg f<\deg_K(a)$. This can
be stated more precisely in the lemma below.

\begin{lm}\label{lm2}\cite[Theorem 2.1(a)]{APZ1}
Fix $(a,\delta)$ a minimal pair.
For any $f\in\overline{K}[X]$ with $\delta_{\overline{\mu}}(f)
<\delta$ we have
\[\init_{\overline{\mu}_{X-a}}(f)=
\init_{\overline{\mu}_{X-a}}(f(a)).\]
Furthermore $\overline{\mu}_{X-a}(f)=\overline{\mu}(f)$.
\end{lm}

We need an auxiliary lemma, in order to prove the second point especially:

\begin{lm}\label{aux}
For any $f\in K[X]$, and $(a,\delta)$ a minimal pair for 
$\overline{\mu}$, such that $\deg f<\deg_K(a)$, then
$\delta_{\overline{\mu}}(f)<\delta$.
\end{lm}

\begin{proof}
Consider any root $b\in\overline{K}$ of $f$. Then 
$\deg_K(b)\leqslant\deg f<\deg_K(a)$, thus by the definition of
minimal pairs, we have $\overline{\mu}(a-b)<\delta$. Now
suppose that $b$ is such that $\delta_{\overline{\mu}}(f)=
\overline{\mu}(X-b)$. By the ultrametric inequality
$\overline{\mu}(X-b)=\overline{\mu}(X-a+(a-b))=
\overline{\mu}(a-b)<\delta$.
\end{proof}

We give a presentation of the proof of \cref{lm2} that differs
from that of \cite{APZ1}, as it shall serve the proof of
\cref{theorem}.

\begin{proof}
Consider the roots of $f$ and write it as
\[f=c\prod_{i=1}^d (X-c_i),\text{ with }c,c_i\in\overline{K}.\]
By hypothesis
\[\forall i,\,\overline{\mu}(X-c_i)\leqslant\delta(f)<\delta=
\delta(Q)=\overline{\mu}(X-a).\]
By the ultrametric property
\begin{align*}
    \overline{\nu}(a-c_i) &=\overline{\mu}((X-c_i)-(X-a))\\
        &=\overline{\mu}(X-c_i)\\
        &<\overline{\mu}(X-a)
\end{align*}
so we naturally have
\[\init_{\overline{\mu}_{X-a}}(X-c_i)=
\init_{\overline{\mu}_{X-a}}(X-a+a-c_i)=
\init_{\overline{\mu}_{X-a}}(a-c_i)\]
thus, by multiplicativity of $\init_{\overline{\mu}_{X-a}}$

\begin{align*}
\init_{\overline{\mu}_{X-a}}(f) &=
\init_{\overline{\mu}_{X-a}}(c)\prod_{i=1}^d
    \init_{\overline{\mu}_{X-a}}(X-c_i)\\
&=\init_{\overline{\mu}_{X-a}}(c)\prod_{i=1}^d
    \init_{\overline{\mu}_{X-a}}(a-c_i)\\
&=\init_{\overline{\mu}_{X-a}}\left(c\prod_{i=1}^d a-c_i\right)\\
&=\init_{\overline{\mu}_{X-a}}(f(a)).
\end{align*}
Furthermore we can extract from this, the following
\begin{align*}
    \overline{\mu}_{X-a}(X-c_i)&=\min\big\{\overline{\mu}(X-a),
        \overline{\nu}(a-c_i)\big\}\\
    &=\overline{\nu}(a-c_i)\\&=\overline{\mu}(X-c_i).
\end{align*}
\end{proof}

\section{Correspondence between Minimal Pairs and Key
Polynomials}

In this section we wish to relate the extensions of 
valuations given by minimal pairs and those given by truncation
by ABKPs. One such stride has been made in the work of
Novacoski. We cite his result in the following theorem.

\begin{thm}\cite[Proposition 3.2]{Nov}
Let $a\in\overline{K}$ be a root of an irreducible polynomial
$Q\in K[X]$ verifying $\delta=\delta(Q)=\overline{\mu}(X-a)$. 
Then
\[Q\text{ is an ABKP for }\mu\iff(a,\delta)
\text{ is a minimal pair for }\overline{\mu}.\]
\end{thm}

\begin{defn}
For any polynomial $f\in\overline{K}[X]$ we shall call
$a\in\overline{K}$ an \emph{optimising root of} $f$, if $a$ is
a root of $f$ and $\overline{\mu}(X-a)=\delta(f)$.
\end{defn}

We will prove that the truncated valuation $\mu_Q$, with $Q$ an
ABKP, comes as a restriction to $K[X]$ of truncated valuations 
on $\overline{K}[X]$, \ie defined via a minimal pair. 

\begin{thm}\label{theorem}
Let $Q\in K[X]$ and let $a\in\overline{K}$ be an optimising
root of $Q$. Then $Q$ is an ABKP, if and only if
$(a,\overline{\mu}(X-a))$ is a minimal pair. Furthermore
$\overline{\mu}_{X-a}$ is an extension of $\mu_Q$
\[\overline{\mu}_{X-a}\big|_{K[X]}=\mu_Q\]
thus inducing a natural injective map of graded algebras
\begin{center}
    \begin{tikzcd}
    \theta\ : \ \text{gr}_{\mu_Q}(K[X])\ar[r, hook]&
        \text{gr}_{\overline{\mu}_{X-a}}(\overline{K}[X]).
    \end{tikzcd}
\end{center}
It sends $\init_{\mu_Q}(f)$ to $\init_{\overline{\mu}_{X-a}}
(f)$, for any $f\in K[X]$.
\end{thm}

This completes the correspondence between the situation over
$K$ and $\overline{K}$. The proof is broken down into three
steps.

\subsection{Step 1.}

We first prove that the valuation given by the minimal pair is
greater in value than the truncated valuation. This is a simple
consequence of the ultrametric inequality.

\begin{lm}\label{ineq}
\[\overline{\mu}_{X-a}\big|_{K[X]}\geqslant\mu_Q.\]
\end{lm}

\begin{proof}
Fix $f\in K[X]$ and set its $Q$-standard decomposition
\[f=\sum_{i}f_i Q^i\]
thus by ultrametric inequality we have
\[\overline{\mu}_{X-a}(f)\geqslant\min_i
\overline{\mu}_{X-a}(f_i Q^i).\]
Now by \cref{lm2}
\[\forall i,\;\overline{\mu}_{X-a}(f_i)=\overline{\mu}(f_i)=
\mu(f_i)\]
and
\[\overline{\mu}_{X-a}(Q)=\overline{\mu}(Q)=\mu(Q).\]
Thus yielding
\[\overline{\mu}_{X-a}(f)\geqslant
\min_i\overline{\mu}_{X-a}(f_i Q^i)=
\min_i\mu(f_i Q^i)=\mu_Q(f).\]
\end{proof}

\subsection{Step 2.}

We build our morphism $\theta$, based on the previous inequality.

\begin{corr}
\cref{ineq} induces a homogeneous ring morphism
\begin{center}
    \begin{tikzcd}
    \text{gr}_{\mu_Q}(K[X])\ar[r,"\theta"]&
        \text{gr}_{\overline{\mu}_{X-a}}(\overline{K}[X])
    \end{tikzcd}
\end{center}
\end{corr}
Indeed, because of \cref{ineq} one has the following inclusions

\begin{align*}
    \mathcal{P}_\gamma(K[X],\mu_Q) &\subseteq
        \mathcal{P}_\gamma(\overline{K}[X],
            \overline{\mu}_{X-a})\\
    \mathcal{P}_\gamma^+(K[X],\mu_Q)&\subseteq
        \mathcal{P}_\gamma^+(\overline{K}[X],
            \overline{\mu}_{X-a})
\end{align*}
for any $\gamma\in\Phi(\mu_Q)\subseteq\Phi(
\overline{\mu}_{X-a})$. The map takes $\init_{\mu_Q}(f)$, which
is an element in degree $\gamma=\mu_Q(f)$, and it sends it to
$f\mod\mathcal{P}_\gamma^+(\overline{K}[X],
\overline{\mu}_{X-a})$. Thus, if
$\overline{\mu}_{X-a}(f)=\mu_Q(f)$, then $\theta$ sends
$\init_{\mu_Q}(f)$ to $\init_{\overline{\mu}_{X-a}}(f)$ and if
$\overline{\mu}_{X-a}(f)>\mu_Q(f)$, $\init_{\mu_Q}(f)$ is sent
to $0$. Furthermore, taking into account its construction,
$\theta$ is homogeneous.

\begin{lm}
We have $\overline{\nu}_{a,\delta}(Q)=\mu(Q)$, thus
$\theta\left(\init_{\mu_Q}Q\right)=
\init_{\overline{\nu}_{a,\delta}}Q\neq0$.
\end{lm}

\begin{proof}
First we have $\delta=\epsilon_\mu(Q)=\max_{i\geqslant1}
\frac{\mu(Q)-\mu(\partial_i Q)}{i}$ thus
\[\mu(Q)=\min_{i\geqslant 1}\{\mu(\partial_i Q)+i\delta\}.\]
Furthermore $0\notin S_{X-a}(Q)$ (because $Q(a)=0$ so
$\overline{\nu}(Q(a))=\infty$), thus
\[\overline{\mu}_{X-a}(Q)=\min_{i\geqslant1}
\big\{\overline{\nu}(\partial_i Q(a))+i\delta\big\}.\]
Thirdly, for any $i\geqslant 1,\;\deg\partial_i Q<\deg Q=
\deg_K(a)$, so that $\delta_{\overline{\mu}}(\partial_i Q)<
\delta$ by \cref{aux}
\[\mu(\partial_i Q)=\overline{\nu}(\partial_i Q(a)).\]
Putting it all together
\[\mu(Q)=\min_{i\geqslant1}\{\mu(\partial_i Q)+i\delta\}         
=\min_{i\geqslant1}\{\overline{\nu}(\partial_i Q(a))+i\delta\}
=\overline{\mu}_{X-a}(Q).\]
\end{proof}

\subsection{Step 3.}

Let us consider the kernel of $\theta$. Considering the 
previous step, one has

\[\text{Ker}(\theta)=\Big\langle\init_{\mu_Q}(I)\Big\rangle,\text{ where }I=
\{f\in K[X]\ |\ \overline{\mu}_{X-a}(f)>\mu_Q(f)\}.\]
Proving \cref{theorem} amounts to showing that $\theta$ is
injective. The main ingredient here is the structure of the graded
algebra. According to \cref{secabkp} these have a polynomial
presentation

\[\text{gr}_{\mu_Q}(K[X]) = G_{<\alpha}[T]\]
where $T=\init_{\mu_Q}(Q),\ \alpha=\deg Q$ and 

\[G_{<\alpha}=\sum_{\deg f<\deg Q}
\text{gr}_{\nu}(K)\init_{\mu_Q}(f).\]
We can make the same statement for $\overline{\mu}_{X-a}$

\[\text{gr}_{\overline{\mu}_{X-a}}(\overline{K}[X])
=(\text{gr}_{\overline{\nu}}(\overline{K}))[\overline{T}]\]
where $\overline{T} = \init_{\overline{\mu}_{X-a}}(X-a)$.
Indeed, by definition $G_{<\deg X-a}$ is generated 
by initial forms of polynomials of degree $<1=\deg(X-a)$, \ie
coefficients in $\overline{K}$. This remark turns out to be crucial
in many ways. One can use \cref{lm2} to prove the following
proposition (which could be seen as a reformulation of \cref{lm2}).

\begin{prop}
The map $\theta$ restricted to $G_{<\alpha}$ induces an injective 
map between $G_{<\alpha}$ and 
$\text{gr}_{\overline{\nu}}(\overline{K})
\subseteq\text{gr}_{\overline{\nu}}(\overline{K})[\overline{T}]$.
\end{prop}

More explicitly, the $\theta$ map will take $\init_{\mu_Q}(f)$ of
$\deg(f)<\deg(Q)$ and send it to $\init_{\overline{\nu}}(f(a))$.
We can now relate our different objects by displaying them in a
commutative diagram

\begin{center}
    \begin{tikzcd}
    G_{<\alpha} \ar[r, hook, "\theta"] \ar[d, hook] &
            \text{gr}_{\overline{\nu}}(\overline{K}) \ar[d, hook]\\
    \text{gr}_{\mu_Q}(K[X])\ar[r,"\theta"] &
            \text{gr}_{\overline{\mu}_{X-a}}(\overline{K}[X])
    \end{tikzcd}
\end{center}
Again by \cref{lm2}, we have that 
$\theta(\init_{\mu_Q}Q)=\init_{\overline{\mu}_{X-a}}Q\neq 0$, because 
$\overline{\mu}_{X-a}(Q)=\mu(Q)=\mu_Q(Q)$. The initial form 
$\init_{\overline{\mu}_{X-a}}$ is multiplicative, so $X-a$ dividing
$Q$ in $\overline{K}[X]$ which implies that
$\init_{\overline{\mu}_{X-a}}(X-a)$ divides
$\init_{\overline{\mu}_{X-a}}(Q)$, \ie $d_{X-a,\overline{\mu}}(Q)>0$.
Hence $\init_{\overline{\mu}_{X-a}}(Q)\notin\text{gr}_{\overline{\nu}}
(\overline{K}).$

We will conclude by using the following basic result on
polynomial rings.

\begin{lm}
Let $\phi\ :\ R\longrightarrow S$ be an injective integral domain
map, that we extend to a map $\Tilde{\phi}\ :\ R[X]\longrightarrow
S[Y]$ which assigns to $X$ a non-constant polynomial $p(Y)\in
S[Y]\setminus S$. Then $\Tilde{\phi}$ is again injective.
\end{lm}

\begin{proof}
Consider a non-zero polynomial $q(X)\in R[X]$. If $q\in R$ then by
assumption $\Tilde{\phi}(q)=\phi(q)\neq 0$. Otherwise suppose
$\deg q(X)\geqslant 1$. Since all the rings we are considering are
integral domains, we get $\deg\Tilde{\phi}(q(X))=\deg q\deg p
\geqslant0\implies\Tilde{\phi}(Q(X))\neq 0$.
\end{proof}

Thus $\theta$ is injective. This concludes our proof.

\subsection{Application to residually transcendental extensions.}

We start by refining \cref{rt}.

\begin{prop}\label{appli}\cite[Theorem 1.3]{Nov}
Our extension $(K[X],\mu)/(K,\nu)$ is residually transcendental if
and only if there is a polynomial $Q$ such that $\mu=\mu_Q$ and
$\epsilon_\mu(Q)$ is torsion over $\Phi(\nu)$\footnote{\ie
$\epsilon_\mu(Q)\in\Q\otimes\Phi(\nu)$}. Furthermore $Q$ can be taken to be
an ABKP.
\end{prop}

\begin{proof}
Choose an extension $\overline{\mu}$ of $\mu$ to $\overline{K}[X]$,
which restricts to $\overline{\nu}$ over $\overline{K}$
($\overline{\nu}$ is thus an extension of $\nu$).
By \cref{rt}, our extension is residually transcendental if and only
if we are given a pair of definition $(a,\delta)\in\overline{K}\times
\Phi({\overline{\nu}})$ for $\overline{\mu}$

\[\overline{\mu}_{X-a}(\sum_i a_i(X-a)^i)=
\min_i\{\overline{\nu}(a_i)+i\delta\},
\text{ where }a_i\in\overline{K}.\]

We can suppose, that this pair is minimal, simply by choosing $a$ 
of minimal degree over $K$, among the pairs defining $\overline{\mu}$.
Thus the minimal polynomial of $Q$ of $a$ is an ABKP and by our
\cref{theorem}, we have that $(K[X],\mu)/(K,\nu)$ is residually
transcendental if and only if
\[\mu=\overline{\mu}\big|_{K[X]}=\overline{\nu}_{a,\delta}\big|_{K[X]}=
\overline{\mu}_{X-a}\big|_{K[X]}=\mu_Q.\]
\end{proof}

\subsection{Application to graded-algebra structure.}

Now recall \cref{sum}. We can now extend the condition under the
sum, defining $G_{<\alpha}$ in the proposition below.

\begin{prop}\label{modif}
Set $f\in K[X],\epsilon(f)<\epsilon(Q)$ where $Q$ is an
ABKP for a valuation $\mu$ over $K[X]$. There exists $g\in
K[X]$ such that $\deg g<\deg Q$ and
$\init_\mu(f)=\init_\mu(g)$. Thus one can re-write the
definition of $G_{<\alpha}$ where $\alpha=\deg Q$:

\[G_{<\alpha}=\sum_{\epsilon(f)<\epsilon(Q)}
\text{gr}_\nu(K)\init_{\mu_Q}(f).\]
\end{prop}

\begin{proof}
We already have

\[G_{<\alpha}\subseteq\sum_{\epsilon(f)<\epsilon(Q)}
\text{gr}_\nu(K)\init_{\mu_Q}(f)\]
since $\deg f<\deg Q$ implies $\epsilon(f)<\epsilon(Q)$ by 
definition of ABKPs. Consider $f$ with $\epsilon(f)<\epsilon(Q)$
and its $Q$-expansion

\[f=\sum_i f_iQ^i,\;\deg f_i<\deg Q.\]
Consider an optimising root of $Q$, $a\in\overline{K}$, thus
providing us with a minimal pair $(a,\delta),\; \delta=
\delta_{\overline{\mu}}(Q)=\epsilon_{\mu}(Q)$.
We thus have the following equalities

\begin{align*}
    \theta\left(\init_{\mu_Q}(f)\right) &= \init_{\overline{\mu}_{X-a}}(f)\\
        &\stackrel{(*)}{=}\init_{\overline{\mu}_{X-a}}(f(a))\\
        &=\init_{\overline{\mu}_{X-a}}(f_0(a))\\
        &\stackrel{(**)}{=}\init_{\overline{\mu}_{X-a}}(f_0)\\
        &=\theta\left(\init_{\mu_Q}(f_0)\right).
\end{align*}
$(*)$ is a direct application of \cref{lm2} and so is $(**)$, since
having $\deg f_0<\deg Q$, implies $\epsilon(f_0)<\epsilon(Q)$.
Thus, since $\theta$ is injective, $\init_{\mu_Q}(f)=
\init_{\mu_Q}(f_0)\in G_{<\alpha}$.
\end{proof}

\subsection{A numerical result.}

We now prove \cref{num}, that for any
$f\in K[X],\:\epsilon_{\mu_Q}(f)\leqslant
\epsilon_{\mu_Q}(Q)=\epsilon_\mu(Q)$. From the proof of
\cref{modif}, one can extract the fact that $S_Q(f)=\{0\}$
whenever  $\epsilon(f)<\epsilon(Q)$. This is proven
in \cite[Prop. 17, Prop. 18]{DMS} and the authors of the
article also announce the proof of its converse in an upcoming
article. We now give a proof of these facts and also prove the
converse of \cite[Prop. 18]{DMS}.

\begin{prop}\label{numBis}
For any polynomial $f\in K[X]$ and any ABKP for a valuation 
$\mu$ over $K[X]$ we have
\[\epsilon_{\mu_Q}(f)\leqslant\epsilon_{\mu_Q}(Q)=
\epsilon_\mu(Q).\]
\end{prop}

\begin{proof}
We will fix $\overline{\mu}$ an extension of $\mu$ to 
$\overline{K}[X]$ and set $a$ an optimising root of $Q$.
By \cref{theorem} we know that $\mu_{X-a}$ is an extension of
$\mu_Q$ and by \cref{thm} we can write
\[\epsilon_{\mu_Q}(f)=\delta_{\mu_{X-a}}(f)=
\max\{\overline{\mu}_{X-a}(X-c);\:c\text{ root of }f\}\]
and $\epsilon_\mu(Q)=\epsilon_{\mu_Q}(Q)=\overline{\mu}(X-a)$.
Thus we simply need to show that for any $c\in\overline{K}$
$\overline{\mu}_{X-a}(X-c)\leqslant\overline{\mu}(X-a)$, which
is straightforward
\[\overline{\mu}_{X-a}(X-c)=\min\{\overline{\mu}(X-a),
\overline{\nu}(a-c)\}\leqslant\overline{\mu}(X-a).\]
\end{proof}

\begin{thm}\label{result}
$\forall f\in K[X],\, \epsilon_{\mu_Q}(f)\leqslant\epsilon_{\mu_Q}(Q)=
\epsilon_\mu(Q)$ and we have equality if and only if there is a $b\geqslant 1$ such
that $\frac{\mu_Q(f)-\mu_Q(\partial_b f)}{b}=\epsilon_\mu(Q)$. In other words
\[S_{Q,\mu}(f)\neq\{0\}\iff\epsilon_{\mu_Q}(f)=\epsilon_\mu(Q).\]
\end{thm}

\begin{proof}
Let us henceforth write $\epsilon_\mu(Q)=\epsilon$. We choose an optimising
root of $Q$, $a\in\overline{K}$ (\ie $\overline{\mu}(X-a)=
\delta_{\overline{\mu}}(Q)$), so that one can apply \cref{theorem}. Since
$\overline{\nu}_{a,\delta}$ is an extension of $\mu_Q$, one can write
\[\epsilon_{\mu_Q}(f)=\epsilon_{\overline{\nu}_{a,\delta}}(f).\]
Let us abbreviate the following initial forms
\begin{align*}
                T &= \init_{\mu_Q}(Q)\\
    \overline{T} &= \init_{\overline{\nu}_{a,\delta}}(X-a).
\end{align*}
We can now express the following equivalences
\begin{align*}
    S_{Q,\mu}(f)\neq\{0\} &\iff d_Q(f)=\deg_{T}\init_{\mu_Q}(f)\geqslant 1\\
        &\iff d_{X-a}(f)=\deg_{\overline{T}}\init_{\overline{\mu}_Q}(f)
                                                    \geqslant 1\\
        &\iff S_{X-a,\overline{\mu}}(f)\neq\{0\}.
\end{align*}
Indeed, recall the properties of the map $\theta$ in the proof of 
\cref{theorem}. It sends initial forms of polynomials $f$ with
$\deg_T\init_{\mu_Q}f=0$ to an initial form with $\deg_{\overline{T}}
\init_{\overline{\mu}_{X-a}}f=0$ and initial forms of polynomials $f$
with $\deg_T\init_{\mu_Q}f\geqslant 1$, to initial forms with
$\deg_{\overline{T}}\init_{\overline{\mu}_{X-a}}f\geqslant 1$.

One can then simply work in $\overline{K}[X]$.
Let $f\in\overline{K}[X]$ such that $S_{X-a,\overline{\mu}}\neq\{0\}$ so
$d=d_{X-a,\overline{\mu}}(f)\geqslant 1$, \ie
\[\overline{\mu}_{X-a}(f)=\overline{\nu}_{a,\delta}(f)=
\overline{\mu}(\partial_d f(a)(X-a)^d).\]
But since we have, by definition of truncation, the inequality
\[\overline{\nu}(\partial_d f(a))\geqslant
\overline{\nu}_{a,\delta}(\partial_d f)\]
which implies
\[\overline{\nu}_{a,\delta}(f)=\overline{\mu}(\partial_d f(a))+
d\epsilon\geqslant\overline{\nu}_{a,\delta}(\partial_d f)+d\epsilon\]
thus
\[\epsilon=\epsilon_{\mu_Q}(Q)\overset{(*)}{\geqslant}
\epsilon_{\mu_Q}(f)=\epsilon_{\overline{\nu}_{a,\delta}}(f)
\geqslant\frac{\overline{\nu}_{a,\delta}(f)-\overline{\nu}_{a,\delta}
(\partial_d f)}{d}\geqslant\epsilon\]
so everything is an equality (the inequality $(*)$ is just due
to \cref{numBis}). Alternatively, one can start with
$\epsilon_Q(f)=\epsilon_{\overline{\nu}_{a,\delta}}(f)<\epsilon$, so that
\[\forall i\geqslant 1,\; \overline{\nu}_{a,\delta}(f)<
\overline{\nu}_{a,\delta}(\partial_i f)+i\epsilon\leqslant
\overline{\nu}(\partial_i f(a))+i\epsilon\]
thus implying that $S_{Q,\mu}(f)=\{0\}$.

\medskip

Let us now show the converse. We suppose, ad absurdum, that we have
a polynomial $f\in\overline{K}[X]$ verifying the two following
statements:
\begin{enumerate}
    \item $\epsilon_{\overline{\nu}_{a,\delta}}(f)=\epsilon$
    \item $S_{X-a,\overline{\mu}}(f)=\{0\}$.
\end{enumerate}
The first statement is equivalent to saying that
\[\forall j\geqslant1,\:\overline{\nu}_{a,\delta}(f)\leqslant
\overline{\nu}_{a,\delta}(\partial_j f)+j\epsilon\tag{1}\]
and that the inequality is an equality for some $j\in\N,\:d\geqslant1$.
We write $d$ the maximal positive integer with this property (so that
$\overline{\nu}_{a,\delta}(f)=\overline{\nu}_{a,\delta}
(\partial_d f)+d\epsilon$). The second statement says that
\[\forall i\geqslant 1,\:\overline{\nu}_{a,\delta}(f)=
\overline{\nu}(f(a))<\overline{\nu}(\partial_i f(a))+i\epsilon.\]
Taking $i=d$ we can use $(1)$ we have
\[\overline{\nu}_{a,\delta}(\partial_d f)+d\epsilon=
\overline{\nu}_{a,\delta}(f)<\overline{\nu}(\partial_d f(a))+d\epsilon\]
or, by simplifying
\[\overline{\nu}_{a,\delta}(\partial_d f)<
\overline{\nu}(\partial_d f(a)).\tag{2}\]
This implies that $S_{X-a,\overline{\mu}}(\partial_d f)\neq\{0\}$. Indeed,
$S_{X-a,\overline{\mu}}(\partial_d f)=\{0\}$ means that
$\overline{\nu}_{a,\delta}(\partial_d f)=\overline{\nu}(\partial_d(f(a))$,
contradicting $(2)$. By the first direction of the statement proven above, we
have that $\epsilon_{\overline{\nu}_{a,\delta}}(\partial_d f)=\epsilon$,
meaning that we have a positive integer $b\geqslant 1$, such that
\[\epsilon=\frac{\overline{\nu}_{a,\delta}(\partial_d f)-
\overline{\nu}_{a,\delta}(\partial_b\partial_d f)}{b}\]
\ie
\[\overline{\nu}_{a,\delta}(\partial_d f)=\overline{\nu}_{a,\delta}
(\partial_b\partial_d f)+b\epsilon.\tag{3}\]
We know that $\partial_b \partial_d=\binom{b+d}{b}\partial_{b+d}$, so by
combining $(3)$ with $(1)$ we have
\[\overline{\nu}_{a,\delta}(f)=\nu\left(\binom{b+d}{b}\right)+
\overline{\nu}_{a,\delta}(\partial_{b+d}f)+(b+d)\epsilon.\]
$\binom{b+d}{b}$ is an integer so $\nu\left(\binom{b+d}{b}\right)\geqslant0$.
But by the inequality in $(1)$ (taken for $j=b+d$) we also have
\[\overline{\nu}_{a,\delta}(f)\leqslant
\overline{\nu}_{a,\delta}(\partial_{b+d} f)+(b+d)\epsilon\]
yielding the reverse inequality $\nu\left(\binom{b+d}{b}\right)\leqslant0$.
In conclusion
\[\overline{\nu}_{a,\delta}(f)=\overline{\nu}_{a,\delta}
(\partial_{b+d}f)+(b+d)\epsilon\]
which contradicts the maximality of $d$.
\end{proof}

\section{Diskoids, towards a geometric interpretation of truncations}

In this last section, we will try and discuss a possible way 
of giving a geometric interpretation for residually transcendental
valuations.

\subsection{Balls.}

We start by studying the simplest possible construction: balls.
We shall fix a valued field $(K,\nu)$ and an extension to its algebraic
closure $(\overline{K},\overline{\nu})$.

\begin{defn}\label{def1}
Fix $a\in\overline{K}$ and $\delta\in\Phi(\overline{\nu})
\cup\{\infty\}$. Then we set
\[D(a,\delta)=\{x\in\overline{K};\:\overline{\nu}(x-a)\geqslant
\delta\}\]
to be the \emph{closed ball centred around $a$ and of radius
$\delta$}.
\end{defn}

\begin{lm}\label{inf}
Let $D=D(a,\delta)$ closed ball. Then for any polynomial
$f\in\overline{K}[X]$, the following minimum is well defined
and can even be explicitly written
\[\min_{x\in D}\overline{\nu}(f(x))=\min_{i\in\N}\{\overline{\nu}
(\partial_i f(a))+i\delta\}.\]
Thus the map $f\longmapsto\min_{x\in D}\overline{\nu}(f(x))$,
that we denote by $\overline{\nu}_D$, is a valuation, which
coincides with the valuation $\overline{\nu}_{a,\delta}$
given by a (not necessarily minimal) defining pair
$(a,\delta)$ .
\end{lm}

\begin{proof}
Write $f$ as
\[f(X)=\sum_{i=0}^n a_i\left(\frac{X-a}{b}\right)^i,
\quad a_i\in\overline{K}\]
where $b\in\overline{K}$ is an element such that 
$\overline{\nu}(b)=\delta$.
Since $x\in D$ if and only if 
$\overline{\nu}\left(\frac{x-a}{b}\right) \geqslant 0$, for
such a $x$ we have by the ultrametric property
\[\overline{\nu}(f(x))\geqslant
\min_{0\leqslant i\leqslant n}\overline{\nu}(a_i).\]
Let us now show that in fact this last value is assumed by
$f$ inside $D$.
We first normalise the expression of $f$ above, \ie we divide
by the $a_i$ of minimal value and we can assume
$\min_{0\leqslant i\leqslant n}\overline{\nu}(a_i)=0$.
Reducing the polynomial's coefficients to the residue field
$\kappa(\overline{\nu})$ we get a non-zero polynomial
$\overline{f}\in\kappa(\overline{\nu})[X]$.
$\kappa(\overline{\nu})$ being algebraically closed, we can
find a value $\overline{y}\in\kappa(\overline{\nu})$ such
that $\overline{f}(\overline{y})=1$.
We can choose $x\in\overline{K}$, such that $\frac{x-a}{b}$ is
in $\Or(\overline{\nu})$ and reduces to $\overline{y}$ in
$\kappa(\overline{\nu})$, thus $\overline{\nu}(f(x))=0$. Necessarily
$x\in D$, since $\overline{\nu}\left(\frac{x-a}{b}\right)
\geqslant0$.
\end{proof}

\begin{rke}
\begin{enumerate}
    \item According to \cref{rt}, the valuation $\nu_D$ is
    residually transcendental.
    \item If $\delta=\infty$, then $D=D(a,\delta)=\{a\}$ so we
    have
    \[\min_{x\in D} \overline{\nu}(f(x)))=\overline{\nu}(f(a)).\]
\end{enumerate}
\end{rke}

We can restate and extend \cref{rt} in the following way.

\begin{thm}\label{rt2}
There is a one-to-one correspondence between residually transcendental
extensions and balls as we have just defined. More precisely, we can
define the following maps, which form a commutative diagram

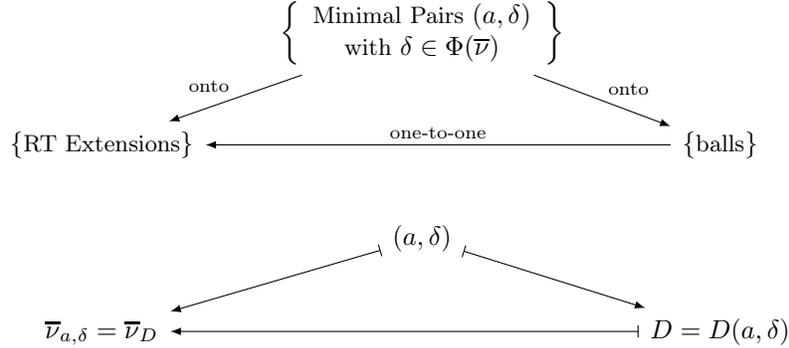
\begin{center}
    \begin{tikzcd}
        &\left\{\begin{array}{c}
            \text{\small Minimal Pairs }(a,\delta)\\
            \text{\small with }\delta\in\Phi(\overline{\nu})
        \end{array}\right\}\ar[rd, "\text{onto}"]\ar[ld,"\text{onto}"']&\\
        \{\text{\small RT Extensions}\}
        &&\{\text{\small balls}\}\ar[ll,"\text{one-to-one}"']\\
        &(a,\delta)\ar[ld,mapsto]\ar[rd,mapsto]&\\
            \overline{\nu}_{a,\delta}=\overline{\nu}_D&&
                D=D(a,\delta)\ar[ll,mapsto]\\
    \end{tikzcd}
\end{center}
\end{thm}

\begin{proof}
We already know that the above maps are well-defined. We also know
that they are onto. Indeed, on one hand, any residually transcendental
extension is given by a minimal pair as we have established in
\cref{rt}. On the other hand, by the ultrametric property we know that
any point inside a ball $D(b,\delta)$ is a centre of the ball, thus we
can choose one of minimal degree, which yields a minimal pair.

Finally the horizontal map is injective. Indeed consider two balls
$D=D(a,\delta),D^\prime(a^\prime,\delta^\prime)$, yielding the same
valuations $\overline{\nu}_D=\overline{\nu}_{D^\prime}$. We can
furthermore suppose that $(a,\delta),(a^\prime,\delta^\prime)$ are
minimal pairs, thus we have
$\overline{\nu}_{a,\delta}=\overline{\nu}_{a^\prime,\delta^\prime}$
and by \cref{pairs}, we have $\delta=\delta^\prime$ and
$\overline{\nu}(a-a^\prime)\geqslant\delta$, which in turn is
equivalent to $D=D^\prime$.
\end{proof}

\begin{rke}
If $\text{rk}(\overline{\mu})=\text{rk}(\overline{\nu})+1$,
then there are elements $\delta\in\Phi(\overline{\mu})\setminus
\Phi(\overline{\nu}),\delta>0$, so that $\forall\gamma\in
\Phi(\overline{\nu}),\,\delta>\gamma$. Thus
\[D(a,\delta)=\{a\}.\]
We thus see that there is no hope for a 1-1 correspondence between
valuations and balls. Indeed $(a,\delta)$ and $(a,\infty)$ give the
same balls, yet their valuations are different:
\[\overline{\nu}_{(a,\delta)}(X-a)=\delta<\infty=
\overline{\nu}_{(a,\infty)}(X-a).\]
The rational rank can increase without the rank increasing, but that
would not affect the correspondence. Consider for instance $K=\C((t))$
equipped with the $t$-multiplicity valuation: $\nu=\text{ord}_t$.
The algebraic closure of $K$ is the field of Puiseux series with
complex coefficients:
$\overline{K}=\bigcup_{n\in\N^*}\C((t^{1/n}))$. We can extend $\nu$
to $\overline{K}$, by extending the $t$-multiplicity. Over 
$\overline{K}[X]$, we can extend $\overline{\nu}$, by considering the
monomial valuation, associating to $X$ the value $\sqrt{2}$. The
resulting valuation $\overline{\mu}$ is still of rank $1$, its
rational rank is $2$, yet there is no problem as the one shown above:
$\forall\delta\in\Phi(\overline{\mu})=\Q(\sqrt{2}),\exists\gamma\in
\Phi(\overline{\nu})=\Q,\,\gamma>\delta$. However, the valuation will
not be RT. This valuation $\mu$ will be out of our correspondence.
In the future we could extend our correspondence to include such
cases as well.
\end{rke}

\subsection{Diskoids.}

Our wish is to generalise this statement up to truncated valuations
$\mu_Q$ which are residually transcendental (\ie such that
$\mu(Q)\in\Q\otimes\Phi(\nu)=\Phi(\overline{\nu})$). More specifically
we wish to find subsets $D\subseteq\overline{K}$ that can verify the
following:

\begin{itemize}[noitemsep]
    \item[(D1)] $\forall f\in K[X],\:\overline{\nu}_D:=
    \min_{x\in D}\overline{\nu}\big(f(x)\big)$ is well-defined.
    \item[(D2)] $\overline{\nu}_D$ is a valuation.
    \item[(D3)] There is a way of associating to $Q$, ABKP for $\mu$ with 
    $\epsilon_\mu(Q)\in\Q\otimes\Phi(\nu)$, a subset $D_Q$ such that 
    $\overline{\nu}_{D_Q}=\mu_Q$.
    \item[(D4)] The mapping $D\longmapsto\overline{\nu}_D$ between such subsets and
    RT extensions over $K[X]$ is bijective.
\end{itemize}

Our candidate was already used by Rüth in his PhD thesis \cite{Ru}
called diskoid. We wish for it to replace the role of balls in the
above correspondence.

\begin{defn}
For any polynomial $f\in\overline{K}[X]$ and any value 
$\rho\in\Phi(\overline{\nu})\cup\{\infty\}$ we define the
\emph{diskoid centred at $f$ of radius $\rho$}
\[\widetilde{D}(f,\rho)=\{x\in\overline{K};\:\overline{\nu}(f(x))
\geqslant\rho\}.\]
\end{defn}

Observe that we then have $\widetilde{D}(X-a,\delta)=D(a,\delta)$. 
Diskoids are a generalisation of balls.
Just like in the previous section, we wish to show that for an
abstract key polynomial $Q$ of value $\mu(Q)\in\Phi(\overline{\nu})$,
the following map
\[f\in K[X]\longmapsto\min_{x\in\widetilde{D}(Q,\mu(Q))}
\overline{\nu}(f(x))\]
is well-defined and corresponds to $\mu_Q$. This is not an easy task,
but we can start by investigating the structure of these diskoids. We
will show that in fact they are disjoint unions of balls.

\begin{rke}
If $\text{rk}(\mu)=\text{rk}(\nu)+1$, \ie when there is an
added convex subgroup when passing from 
$(\overline{K},\overline{\nu})$ to
$(\overline{K}[X],\overline{\mu})$, we have key
polynomials $Q$ such that
$\mu(Q)\in\Phi(\overline{\mu})\setminus\Phi(\overline{\nu}),\:
\mu(Q)>0$, so that
$\forall\gamma\in\Phi(\overline{\nu}),\,\gamma<\mu(Q)$.
In this case diskoids are reduced to finite sets of points:
for any $x\in\overline{K},\,\overline{\nu}(f(x))\in
\Phi(\overline{\nu})$, so $\overline{\nu}(f(x))<\mu(Q)$,
unless $x$ is a root of $f$, in which case
$\overline{\nu}(f(x))=\infty$.
Thus there is no hope for a 1-1 correspondence. For instance,
if we consider only valuations $\mu$ of rank $1$, such problem will
not arise. Indeed, in that case $\Phi(\overline{\mu})$ can be
embedded as an ordered group inside $\R$.

If $\mu$ is residually transcendental, such problems do not 
arise either, since in this case 
$\Phi(\overline{\mu})=\Phi(\overline{\nu})$ and
$\text{rk}(\mu)=\text{rk}(\nu)$.
\end{rke}

We need the following lemma.

\begin{lm}\label{mufactor}
Set $f\in\overline{K}[X],\; a\in\overline{K}$ a root of $f$ and 
$\rho\in\Phi(\overline{\nu})$. Then the quantity
\[\epsilon(a;f,\rho)=\min\{\lambda\in\Phi(\overline{\nu});\:
D(a,\lambda)\subseteq\widetilde{D}(f,\rho)\}\]
is well-defined and can even be explicitly written
\[\epsilon(a;f,\rho)=\max_{i\in\N^*}\ \frac{\rho-\overline{\nu}
(\partial_i f(a))}{i}.\]
\end{lm}

\begin{proof}
Let us write $f$ as

\[f(X)=\sum_{i=1}^n a_i(X-a)^i,\quad a_i\in \overline{K}\]
so that $a_i=\partial_if(a)$. We remark that $a_0=0$ since $f(a)=0$,
so considering \Cref{inf}, for any $\lambda\in\Phi(\overline{\nu})$,
$\min_{x\in D(a,\lambda)}\overline{\nu}(f(x))=\min_{i\geqslant 1}
\{\overline{\nu}(a_i)+i\lambda\}$. Thus we can write the following
equivalences

\begin{align*}
    D(a,\lambda)\subseteq\widetilde{D}(f,\rho)
        &\iff \forall x\in D(a,\lambda),\ 
            \overline{\nu}(f(x))\geqslant\rho\\
        &\iff\min_{x\in D(a,\lambda)} \overline{\nu}(f(x))\geqslant\rho\\
        &\iff\min_{1\leqslant i\leqslant n}
            \{\overline{\nu}(a_i)+i\lambda\}\geqslant\rho\\
        &\iff\forall i=1,\ldots,n,\:
            \overline{\nu}(a_i)+i\lambda\geqslant\rho\\
        &\iff\forall i=1,\ldots,n,\:
            \lambda\geqslant\frac{\rho-\overline{\nu}(a_i)}{i}\\
        &\iff \lambda\geqslant\max_{1\leqslant i\leqslant n}
            \frac{\rho-\overline{\nu}(a_i)}{i}.
\end{align*}
Thus we can set the following expression
\[\epsilon(a;f,\rho)=\max_{1\leqslant i\leqslant n}
\frac{\rho-\overline{\nu}(a_i)}{i}.\]
\end{proof}

We can now show how diskoids decompose.

\begin{lm}\label{dec}
For any $f\in\overline{K}[X],\, \rho\in\Phi(\overline{\nu}),\:
\widetilde{D}(f,\rho)$ is a finite union of balls. They are centred
around the roots of $f$:
\[\widetilde{D}(f,\rho)=\bigcup_{f(c)=0}D(a,\epsilon(c;f,\rho)).\]
\end{lm}

\begin{proof}
Let $c_1,\ldots,c_n$ be the possibly repeated roots of $f$ in
$\overline{K}$. By \cref{mufactor} we can assign to each $c_i$
the value $\epsilon_i=\epsilon(c_i;f,\rho)\in
\Phi(\overline{\nu})$ such that $D(c_i,\epsilon_i)\subseteq
\widetilde{D}(f,\rho)$, with $\epsilon_i$ minimal for this
property. We now show that $\widetilde{D}(f,\rho)\subseteq
\cup_iD(c_i,\epsilon_i)$, the other inclusion being clear.
Let $x\in\widetilde{D}(f,\rho)$
\[\sum_{i=1}^n \overline{\nu}(x-c_i)=\overline{\nu}(f(x))\geqslant
\rho.\]
We can rearrange the $c_i$ so that
\[\overline{\nu}(x-c_1)\geqslant\ldots
\geqslant\overline{\nu}(x-c_n).\]
In other words we have relabelled the roots to make $c_1$ the 
closest root to $x$. Now we will show that
$D(c_1,\overline{\nu}(x-c_1))\subseteq
D(f,\rho)$ so that $\epsilon_i\leqslant\overline{\nu}(x-c_1)$.
For $y\in D(c_1,\overline{\nu}(x-c_1))$ we have

\begin{align*}
    \overline{\nu}(f(y)) &= \sum_{i=1}^n\overline{\nu}(y-c_i)\\
        &\geqslant\overline{\nu}(x-c_1)+\sum_{i=2}^n
            \overline{\nu}(\underbrace{y-c_1}_{\geqslant
            \overline{\nu}(c_1-x)}+\underbrace{c_1-x}_{\geqslant
            \overline{\nu}(x-c_i)}+x-c_i)\\
        &\geqslant\overline{\nu}(x-c_1)+\sum_{i=2}^n \overline{\nu}(x-c_i)\\
        &\geqslant\rho.
\end{align*}
\end{proof}

\begin{rke}
If $(a,\delta)$ is a minimal pair associated to an abstract key polynomial
$Q$, then
\[\epsilon(a;Q,\mu(Q))=\delta.\]
Indeed, since for any $i\geqslant 1,\, \deg\partial_iQ<\deg Q$, we get by 
\cref{result} $\mu(\partial_i Q)=\mu_Q(\partial_iQ)=
\overline{\nu}_{a,\delta}(\partial_iQ)=\overline{\nu}(\partial_i Q(a))$,
since $\epsilon_{\mu_Q}(\partial_iQ)=\epsilon_\mu(\partial_iQ)<
\epsilon_\mu(Q)$.
\end{rke}

\begin{defn}
For any polynomial $f,g\in\overline{K}[X]$ and value
$\rho\in\Phi(\overline{\nu})$, the following value
\[\min_{x\in\widetilde{D}(f,\rho)}\overline{\nu}(g(x))\]
is well-defined, according to \cref{mufactor} and \cref{inf} (Diskoids
satisfy (D1)). We can thus define the following map
\begin{align*}
    \overline{\nu}_{\widetilde{D}(f,\rho)}:K[X]&
        \longrightarrow\Phi(\overline{\nu})\cup\{\infty\}\\
    g&\longmapsto\min_{x\in\widetilde{D}(f,\rho)}\overline{\nu}(g(x)).
\end{align*}
\end{defn}

\begin{prop}
The map $\overline{\nu}_{\widetilde{D}(f,\rho)}$ is ultrametric
(verifies (V2)).
\end{prop}

\begin{proof}
Take $g,h\in\overline{K}[X]$ and set $x\in\widetilde{D}(f,\rho)$ such
that $\overline{\nu}_{\widetilde{D}(f,\rho)}(g+h)=
\overline{\nu}(g(x)+h(x))$. We have the following
\begin{align*}
    \overline{\nu}_{\widetilde{D}(f,\rho)}(g+h) &= \overline{\nu}(g(x)+h(x))\\
    & \geqslant\min\{\overline{\nu}(g(x)),\overline{\nu}(h(x))\}\\
    &\geqslant\min\{\overline{\nu}_{\widetilde{D}(f,\rho)}(g),
        \overline{\nu}_{\widetilde{D}(f,\rho)}(h)\}.
\end{align*}
\end{proof}

\begin{rke}
In a similar way one can prove that
\[\forall g,h\in\overline{K}[X],\:
\overline{\nu}_{\widetilde{D}(f,\rho)}(gh)
\geqslant\overline{\nu}_{\widetilde{D}(f,\rho)}(g)+
\overline{\nu}_{\widetilde{D}(f,\rho)}(h).\]
However $\overline{\nu}_{\widetilde{D}(f,\rho)}$ may fail to be
multiplicative (the condition (V1) for valuations) so diskoids fail to
satisfy (D2). We give a simple counter-example: set $a\in K$ such that
$\nu(a)<0$. Define $D_a=D(a,0),\; D_0=D(0,0)$ and $D=D_a\cup D_0$.
It can be realised as a diskoid of $f=X(X-a)$
\[\widetilde{D}(X(X-a),\nu(a))=D(a,0)\cup D(0,0)=D\]
since
\[\epsilon(0;f,\nu(a))=\epsilon(a;f,\nu(a))=
\max\left\{\frac{\nu(a)}{2},\nu(a)-\nu(a)\right\}=0.\]
Indeed $\partial_1 X(X-a)=2X-a$ and $\partial_2 X(X-a)=1$. Then we
have
\[\begin{array}{ll}
    \overline{\nu}_{D_0}(X)=0 & \overline{\nu}_{D_a}(X)=\nu(a) \\
    \overline{\nu}_{D_0}(X-a)=\nu(a) & \overline{\nu}_{D_a}(X-a)=0
\end{array}\]
so clearly
\[\overline{\nu}_{\widetilde{D}(f,\overline{\nu}(a))}(X)=
\overline{\nu}_{\widetilde{D}(f,\overline{\nu}(a))}(X-a)=\nu(a).\]
Furthermore $\overline{\nu}_{\widetilde{D}(f,\overline{\nu}(a))}
(X(X-a))=\nu(a)$ but
\[\overline{\nu}_{\widetilde{D}(f,\overline{\nu}(a))}(X)+
\overline{\nu}_{\widetilde{D}(f,\overline{\nu}(a))}(X-a)=
2\nu(a)<\nu(a)=\overline{\nu}_{\widetilde{D}(f,\nu(a))}(X(X-a))\]
hence, $\overline{\nu}_{\widetilde{D}(f,\overline{\nu}(a))}$ is not a
valuation.
\end{rke}

One can hope that if we restrict to studying only abstract key
polynomials $Q$, with $\mu(Q)\in\Phi(\overline{\nu})$, the diskoids
$\widetilde{D}(Q,\mu(Q))$ yield true valuations. One could even
generalise \cref{rt2} and interpret the decomposition of a diskoid
into balls at the level of valuations. Let us state these as
conjectures.

\begin{con}

\begin{itemize}[noitemsep]
    \item[(C1)] If $Q$ is an ABKP for $\mu$, then
    $\overline{\nu}_{\widetilde{D}(Q,\mu(Q))}$ is a valuation.
    \item[(C2)] If $\mu$ is an RT extension, there is a unique diskoid $D$ such that
    $\mu=\overline{\nu}_D$.
\end{itemize}
\end{con}

In order to explore these questions, we will recall some basic facts
about henselization.

\subsection{Recall on henselian valued fields and henselization.}

In this section we recall some facts about the interaction between
extensions of valued fields and the Galois theory of the corresponding
field extension. We start with the following crucial notion.

\begin{defn}
A valued field $(K,\nu)$ is called \emph{henselian} if there
is only one extension of $\nu$ to the algebraic closure of 
$K$.
\end{defn}

Several ways to characterise henselian fields exist and relate to
Hensel's lemma or the implicit function theorem for valued fields (one
can consult \cite{Ku} for a deeper understanding).

\begin{ex}
Henselian fields include discrete valued fields of rank $1$ which are
complete, such as the field of $p$-adic numbers $\Q_p$ equipped with
its natural $p$-adic absolute value, or the formal power series
$k((t))$ where $k$ is any field, which we equip with the natural
$t$-adic valuation.
\end{ex}

Not every field is henselian, but one can find a smallest henselian
extension of a valued pair $(K,\nu)$ which is henselian. It is called
a henselization of $(K,\nu)$.

\begin{defn}
An extension $(\Tilde{K},\Tilde{\nu})$ of $(K,\nu)$ is called a 
\emph{henselization} of $(K,\nu)$ if it is henselian and if for every
henselian valued field $(E,\zeta)$ and every embedding $\lambda:
(K,\nu)\hookrightarrow(E,\zeta)$ there exists a unique embedding
$\Tilde{\lambda} : (\Tilde{K},\Tilde{\nu})\hookrightarrow(E,\zeta)$
extending $\lambda$.
\end{defn}

Henselizations exist and can be constructed in the following 
way. Choose a separable closure $K^{sep}$ of $K$ and an
extension $\nu^s$ of $\nu$ to $K^{sep}$ (or $\overline{K}$).
Write $G_K=\text{Gal}(K^{sep}/K)=\text{Aut}_K(\overline{K})$
and
\[G^h=\{\sigma\in G_K\ |\ \nu^s\circ\sigma=\nu^s\}\]
the decomposition group of $\nu^s$. The decomposition field of $\nu^s$
which is by definition
\[K^{h(\nu^s)}=\{x\in K^{sep}\ |\ \sigma(x)=x,\:
\forall\sigma\in G^h\}\]
is a henselization of $(K,\nu)$. Any other choice for the extension
$\nu^s$ will give another henselization. All henselizations are
isomorphic up to unique isomorphism so we will talk about \emph{the}
henselization of $(K,\nu)$ and write it $(K^h,\nu^h)$.

In our situation we are interested in certain types of extensions.
Once we have our truncated valuation $\mu_Q$, we wish to enumerate the
valuations that extend both $\mu_Q$ and $\overline{\nu}$, to
$\overline{K}[X]$. We use the following result concerning the
permutation of valuations by automorphisms of extensions. Consider
$(K,\nu)$ a valued field and $L/K$ a field extension. We will write
$\mathcal{E}(L,\nu)$ the set of valuations on $L$ extending $\nu$.

\begin{prop}\label{perm}
Let $L/K$ be a normal field extension and $\nu$ a valuation of $K$.
$\text{Aut}_K(L)$ acts \emph{transitively} on the set of extensions of
$\nu$ to $L$ as follows
\begin{align*}
    \text{Aut}_K(L)\times\mathcal{E}(L,\nu) &\longrightarrow \mathcal{E}(L,\nu)\\
    (\sigma,\mu) &\longmapsto \mu\circ\sigma^{-1}=\mu^\sigma.
\end{align*}
\end{prop}

Once we have a common extension $\overline{\mu}$ of our $\mu_Q$ and 
$\overline{\nu}$, so for instance $\overline{\nu}_{a,\delta}$ where
$(a,\delta)$ is a minimal pair associated to $Q$, then we can let
$G_K=\text{Aut}_K(\overline{K})$ act on $\overline{\mu}$ so that we
obtain all the extensions of $\mu_Q$. Indeed $\overline{K}(X)/K(X)$ is
a normal extension with same Galois group as $\overline{K}/K$. So, if
$\sigma\in G_K$, we have the following equivalences
\begin{align*}
    \overline{\mu}\circ\sigma^{-1}\in\mathcal{E}(\overline{K}[X],
        \overline{\nu})
    &\iff\overline{\mu}\circ\sigma^{-1}\big|_{\overline{K}}=
        \overline{\nu}\\
    &\iff\overline{\mu}\big|_{\overline{K}}\circ\sigma^{-1}
        \big|_{\overline{K}}=\overline{\nu}\quad
        \text{ since }\sigma^{-1}(\overline{K})\subset\overline{K}\\
    &\iff\overline{\nu}\circ\sigma^{-1}\big|_{\overline{K}}=
        \overline{\nu}\\
    &\iff \sigma^{-1}\big|_{\overline{K}}\in G^h.
\end{align*}
Thus, if we identify the automorphism groups of $\overline{K}(X)/K(X)$
and $\overline{K}/K$ we can state the proposition below.

\begin{prop}
Take an abstract key polynomial $Q$ and $(a,\delta)$ an associated
minimal pair. We thus have

\begin{align*}
\mathcal{E}(\overline{K}[X],\overline{\nu})\cap
\mathcal{E}(\overline{K}[X],\mu_Q)
&=\{\overline{\nu}_{a,\delta}\circ\sigma^{-1};\:\sigma\in G^h\}\\
&=\{\overline{\nu}_{\sigma(a),\delta};\:\sigma\in G^h\}.
\end{align*}
\end{prop}

\begin{proof}
It remains to prove the following identity
\[\overline{\nu}_{a,\delta}\circ\sigma^{-1}=
\overline{\nu}_{\sigma(a),\delta}\]
for any $\sigma\in G^h$. For $f\in\overline{K}[X]$ decompose
it as $f=\sum_{i\geqslant 0}a_i(X-\sigma(a))^i$

\begin{align*}
\overline{\nu}_{a,\delta}\circ\sigma^{-1}(f) &=
    \overline{\nu}_{a,\delta}\left(\sum_{i\geqslant 0}
        \sigma^{-1}(a_i)(X-a)^i\right)\\
    &=\min_{i\geqslant 0}
        \{\overline{\nu}(\sigma^{-1}(a_i))+i\delta\}\\
    &=\min_{i\geqslant 0}\{\overline{\nu}(a_i)+i\delta\} \\
    &=\overline{\nu}_{\sigma(a),\delta}(f).
\end{align*}
\end{proof}

\begin{rke}
Just as $G^h$ acts on valuations, it also acts on balls. Set 
$\sigma\in G^h$ and $(a,\delta)\in\overline{K}\times\Phi(\overline{\nu})$. Then

\[\sigma\left(D(a,\delta)\right)=D(\sigma(a),\delta).\]
Indeed, for any $x\in\overline{K}$, one has $\overline{\nu}(x-a)=\overline{\nu}
(\sigma(x)-\sigma(a))$.
\end{rke}

\subsection{Return to diskoids.}

We can use the action of $G^h$ to make the structure of diskoids easier
to grasp and make more sense of the maps $\overline{\nu}_{\widetilde{D}(f,\mu(f))}$.

\begin{prop}
Consider a polynomial $f\in K^h[X]$, such that $G^h$ acts transitively
on its roots and $\overline{\mu}(f)\in\Phi(\overline{\nu})$. If we set
$a\in\overline{K}$ an optimising root of $f$ (\ie $\overline{\mu}(X-a)=
\delta(f)$), then we have

\[\widetilde{D}(f,\overline{\mu}(f))=
\bigcup_{\sigma\in G^{h(\overline{\nu})}}D(\sigma(a),\delta(f)).\]
Thus $\overline{\nu}_{\widetilde{D}(f,\overline{\mu}(f))}$ is in fact a
valuation on $K[X]$ (or even $K^h[X]$) and  its extensions to
$\overline{K}[X]$ are $\overline{\nu}_{D(\sigma(a),\delta(f))},\,
\sigma\in G^h$.
\end{prop}

\begin{proof}
Considering \cref{dec}, the union is well-indexed, seeing as the roots of
$f$ will be $\sigma(a),\ \sigma\in G^h$. We already know
$\overline{\nu}_{D(\sigma(a),\delta(f))}$ is  well-defined valuation and
equal to the valuation given by the defining pair $(a,\delta(f))$. Consider
any $g\in K^h[X]$, we will show that in fact
$\overline{\nu}_{D(\sigma(a),\delta(f))}(g)=
\overline{\nu}_{D(a,\delta(f))}(g)$ this value thus being also equal to
$\overline{\nu}_{D(f,\overline{\mu}(f))}(g)$. Indeed, for any $\sigma\in
G^h$, $\sigma$ will not change the coefficients of $g$, thus we can safely
write

\begin{align*}
    \overline{\nu}_{D(\sigma(a),\delta(f))}(g) &= 
                \overline{\nu}_{\sigma(a),\delta(f)}(g)\\
    &= \overline{\nu}_{a,\delta(f)}\circ\sigma^{-1}(g)\\
    &= \overline{\nu}_{a,\delta(f)}(g)\\
    &= \overline{\nu}_{D(a,\delta(f))}(g).
\end{align*}
Since $\overline{\nu}_{D(a,\delta(f))}$ is a valuation over
$\overline{K}[X]$, extending $\overline{\nu}_{\widetilde{D}
(f,\overline{\mu}(f))}$, by \cref{perm}, we know that all the other
extensions are $\overline{\nu}_{a,\delta(f)}\circ\sigma^{-1}=
\overline{\nu}_{D(\sigma(a),\delta(f))}$.
\end{proof}

For instance if the polynomial $f$ in the preceding theorem is irreducible,
then it has its roots permuted transitively.

\begin{corr}
Let $Q$ be an abstract key polynomial with $\epsilon_\mu(Q)\in\Phi(\overline{\nu})$.
If $Q$ is irreducible over $K^h$, then $\overline{\nu}_{D(Q,\mu(Q))}$ is a
valuation and in fact

\[\mu_Q=\overline{\nu}_{\widetilde{D}(Q,\mu(Q))}.\]
\end{corr}

\begin{proof}
By \cref{theorem}, if we assign to $Q$ one of its minimal pairs $(a,\delta)$

\begin{align*}
    \mu_Q &= \overline{\nu}_{a,\delta}\big|_{K[X]}\\
    &= \overline{\nu}_{D(a,\delta)}\big|_{K[X]}\\
    &= \overline{\nu}_{\widetilde{D}(Q,\mu(Q))}.
\end{align*}
\end{proof}

We thus formulate the following concept.

\begin{defn}
We say that a polynomial $f\in K[X]$ is \emph{analytically irreducible} if
it is irreducible over $K^h$
\end{defn}

If $Q$ is analytically irreducible, we can then extend the correspondence of
\cref{rt2}, into the one below.

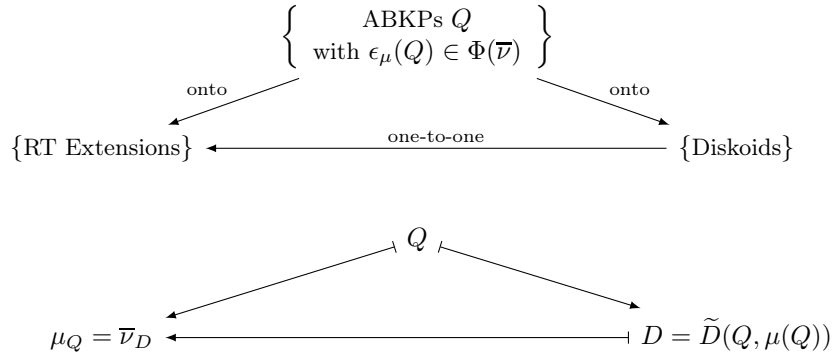
\begin{center}
    \begin{tikzcd}
        &\left\{\begin{array}{c}
            \text{\small ABKPs }Q\\
            \text{\small with }\epsilon_{\mu}(Q)\in\Phi(\overline{\nu})
        \end{array}\right\}\ar[rd, "\text{onto}"]\ar[ld,"\text{onto}"']&\\
        \{\text{\small RT Extensions}\}
        &&\{\text{\small Diskoids}\}\ar[ll,"\text{one-to-one}"']\\
        &Q\ar[ld,mapsto]\ar[rd,mapsto]&\\
            \mu_Q=\overline{\nu}_D&&
                D=\widetilde{D}(Q,\mu(Q))\ar[ll,mapsto]\\
    \end{tikzcd}
\end{center}

If the base field $(K,\nu)$ is itself henselian, the answer is obviously 
affirmative (since ABKPs are irreducible over $K$ to begin with), thus the
correspondence is true in this case.

\subsection{The rank one case.}

This section is dedicated to showing the following theorem.

\begin{thm}
Suppose $\nu$ is of rank $1$ and $Q$ is an ABKP for a valuation $\mu$ of $K[X]$
such that $\mu_Q/\nu$ is an RT extension,\ie $\epsilon_\mu(Q)\in\Q\otimes\Phi(\nu)$.
Then $Q$ is analytically irreducible.
\end{thm}

The reasoning can be done by reductio ad absurdum. Suppose on the contrary
that $Q$ is reducible in $K^h[X]$. We can write $Q=PR,\ P,R\in K^h[X]$ with
$P$ irreducible over $K^h$ and $\epsilon_\mu(P)=\epsilon_\mu(Q)$.
Indeed consider $a$ an optimising root of $Q$ and suppose $P$ is its minimal
polynomial over $K^h$. Then $P|Q$ and since $Q$ is not irreducible, we 
get that $P\neq Q$.

In rank one, the completion satisfies Hensel's lemma (\cf \cite[Theorem 1.3.1]{EP}),
thus it is henselian. This shows that the henselization can be embedded in the
completion, hence $K$ is dense in $K^h$, since $K$ is dense in its completion.
This means that
\[\forall z\in K^h,\:\forall\epsilon\in\Phi(\nu),\:\exists z^*\in K,\:
\nu^h(z-z^*)>\epsilon.\]

We now use the continuity of roots of a monic polynomial of constant degree
relative to its coefficients. We use as a reference \cite[Theorem 2]{Br},
but more specifically the remark that follows Theorem 2 of the
aforementioned paper.

\begin{thm}
Let $f$ and $f^*$ be monic polynomials of common degree $n>1$ with
coefficients in an algebraically closed, valued field $(F,v)$. We suppose
that these coefficients are of positive value. We may then write 
\begin{align*}
    f&=\prod_{k=1}^n(X-\alpha_k)\\
    f^*&=\prod_{k=1}^n(X-\alpha_k^*)
\end{align*}
such that $\nu(\alpha-\alpha^*)\geqslant V(f-f^*)/n$, where
\[V\left(\sum_l a_lX^l\right)=\min_l v(a_l).\]
\end{thm}

Let us come back to our proof. We can arbitrarily approach each
coefficient of $P$ by elements in $K$, thus one can find a monic
$P^*\in K[X]$ of same degree, such that $V(P-P^*)>n\epsilon_\mu(P)$.
We decompose them in linear factors
\begin{align*}
    P&=\prod_{k=1}^r(X-a_k)\\
    P^*&=\prod_{k=1}^r(X-a_k^*)
\end{align*}

so that $\forall k,\:\overline{\nu}(a_k-a_k^*)>\epsilon_\mu(P)$.
Thus $\forall k,\:\overline{\mu}(X-\alpha_k)=\overline{\mu}(X-\alpha_k^*)$,
so the polynomial $P^*\in K^h$, is such that
$\epsilon(P^*)=\epsilon(P)=\epsilon(Q)$ but $\deg(P^*)=\deg(P)<\deg(Q)$.
This contradicts the fact that $Q$ is an ABKP.

\medskip

In higher rank cases, this proof breaks down, since $K^h$ is not longer in
the completion of $K$. In a private communication with F.-V. Kuhlmann, he
expressed his pessimism about analytic irreducibility of ABKPs.

\end{document}